\documentclass[12pt]{paper}
\usepackage{amscd}
\usepackage{amssymb}
\usepackage{a4wide}
\usepackage{amstext}
\usepackage{amsthm}
\usepackage{xcolor}
\usepackage{amsfonts, amsmath, comment}

\usepackage{xcolor}
\usepackage{mathabx}

\usepackage[a4paper, hmargin=2.85cm, vmargin={3.95cm, 3.95cm}]{geometry}

\numberwithin{equation}{subsection}

\newcommand{\Z}{\mathbb{Z}}

\newcommand{\N}{\mathbb{N}}
\newcommand{\R}{\mathbb{R}}

\newcommand{\Ho}{\mathcal{H}}
\newcommand{\So}{\mathcal{S}}

\newcommand{\summ}{\sum\limits}

\newcommand{\eps}{\varepsilon}

\newcommand{\pho}{\rho}

\newcommand{\la}{\lambda}
\newcommand{\Lam}{\Lambda}

\DeclareMathOperator{\supp}{supp} 
\DeclareMathOperator{\sign}{sign}

\def\done{{1\hskip-2.5pt{\rm l}}}

\renewcommand{\phi}{\varphi}

\newtheorem{Thm}{Theorem}
\newtheorem{theorem}[Thm]{Theorem}
\newtheorem{lemma}{Lemma}
\newtheorem{claim}{Claim}

\newtheorem{corollary}{Corollary}

\newtheorem{definition}{Definition}
\begin{document}

\title{Fourier uniqueness and non-uniqueness pairs}
\author{Aleksei Kulikov, Fedor Nazarov, Mikhail Sodin}

\maketitle
\noindent{\sf \today}

\begin{abstract}
Motivated by recent works by Radchenko and Viazovska and by Ramos and Sousa,
we find sufficient conditions for a pair of discrete subsets of the real line
to be a uniqueness or a non-uniqueness pair for the Fourier transform. These conditions are
close to each other.
The uniqueness result can be upgraded to an interpolation formula, which in turn produces an abundance of discrete measures with discrete Fourier transform.
\end{abstract}

\bigskip

\hfill{\sf To Vladimir Marchenko on the occasion of his 101st birthday}

\section{Introduction and  main results}

\subsection{Uniqueness and non-uniqueness pairs}
Let $\mathsf H$ be a space of functions on $\mathbb R$.
For instance, $\mathsf H$ can be $L^2(\mathbb R)$ or the Schwartz space
$\mathcal S=\mathcal S(\mathbb R)$. Later we will also consider other spaces $\mathsf H$.
We define the Fourier transform by
\[
\widehat{f}(\xi) = \int_{\mathbb R} f(x) e^{-2\pi {\rm i} \xi  x}\, {\rm d} x\,.
\]
Then its inverse is given by
\[
f(x) = \int_{\mathbb R} \widehat{f}(\xi) e^{2\pi {\rm i} \xi  x}\, {\rm d} \xi\,.
\]
\begin{definition}[{\sf uniqueness and non-uniqueness pairs}]\label{def-pairs}
We call a pair $(\Lambda, M)$ of sets $\Lambda, M \subset \mathbb R$
a uniqueness pair {\sf (UP)} for $\mathsf H$
if
\[
f\in\mathsf H, \ f|_\Lambda = 0, \ \widehat{f}|_M = 0 \ \Longrightarrow f = 0.
\]
Otherwise, $(\Lambda, M)$ is a non-uniqueness pair {\sf (NUP)}
for $\mathsf H$.
\end{definition}

\noindent
This definition is invariant with respect to the transformations $(\Lambda, M) \mapsto (t\Lambda, t^{-1}M)$, $t\ne 0$, and $(\Lambda, M) \mapsto (\Lambda + r, M + s)$, $r, s \in \mathbb R$,
provided that, for every $f\in\mathsf H$ and every $t\ne 0$, $r, s \in \mathbb R$,
the space $\mathsf H$ also contains the functions $x\mapsto f(tx)$,
$x\mapsto f(x+r)$, and $x\mapsto e^{2\pi {\rm i} sx}f(x)$.

\smallskip
The most classical examples of the uniqueness and non-uniqueness pairs
(for instance, in the space $\mathcal{S}$)  are

\medskip\noindent ({\sf UP}) \quad
$ \Lambda = \mathbb Z$, $M = \mathbb R\setminus (-\tfrac12, \tfrac12)$,

\medskip\noindent ({\sf NUP}) \quad
$\Lambda = c\Z$, $M = \mathbb R\setminus (-\tfrac12, \tfrac12)$ with any $c>1$.

\medskip\noindent
The condition $M = \mathbb R\setminus (-\tfrac12, \tfrac12)$  implies that the spectrum of
$f$ is contained in the interval
$[-\tfrac12, \tfrac12]$, that is, $f$ is an entire function of exponential type at most $\pi$.
The famous Beurling--Malliavin theorem~\cite{BM} can be viewed as a description of the sets $\Lambda$ such that $(\Lambda, \mathbb R\setminus (-a, a)$ is a UP for each $a<\tfrac12$ and a NUP for each $a>\tfrac12$.

\smallskip
Another classical example of a uniqueness pair in $L^2(\mathbb R)$ is due to Benedicks~\cite{Benedicks} and Amrein and Berthier~\cite{AB}:

\medskip\noindent ({\sf UP}) \quad
$ \Lambda$, $M$  have complements of finite Lebesgue measure.

\medskip\noindent
Later, the second-named author~\cite[Section~2.5]{Nazarov}
found ``a one-sided version'' of this result:

\medskip\noindent ({\sf UP}) \quad
$ \Lambda \cap \mathbb R_-$, $M\cap \R_-$
have complements of finite Lebesgue measure.

\medskip\noindent
In all these examples, at least one of the sets $\Lambda$, $M$ in the uniqueness pair was not discrete.

The situation changes if we assume that the functions in $\mathsf H$,
as well as their Fourier transforms, are well-defined pointwise.
In 2017 Radchenko and Viazovska~\cite{Rad-Via} found a surprising example of a
uniqueness pair for the set of {\em even} functions in $\mathcal S$:

\medskip\noindent ({\sf UP}) \quad
$\Lambda = M = \{ 0, \pm \sqrt 1, \pm \sqrt 2, \pm \sqrt 3, \ldots \, \}$.

\medskip\noindent
This follows from their remarkable interpolation formula, which expresses any
function $f\in\mathcal S$ by the values of $f$ and $\widehat f$ on this set and two more values
$f'(0)$ and $\widehat{f}\,'(0)$. Radchenko and Viazovska used an arithmetic approach based on the theory of modular forms. This makes their construction quite rigid.

\smallskip
Our work was inspired by a work of Ramos and Sousa~\cite{Ram-Sousa1}, in which
they developed a purely analytic approach to the construction of discrete uniqueness pairs.
Starting with a function $f\in\mathcal S$ such that $f|_\Lambda$ and $\widehat{f}|_M$ vanish,
their idea was, first, to estimate the size of the derivatives $f^{(j)}$ and $\widehat{f}^{(j)}$, $j=1, 2, \ldots\, $, to show that $f$ and $\widehat f$ have an analytic continuation to $\mathbb C$ with certain bounds, and then, using the entire function theory, that $f$ is a zero function. Their approach used a repeated application of the Rolle's theorem, and required unnecessarily restrictive assumptions on the denseness of $\Lambda$ and $M$. For instance, they were able to show that the pair of the sets
$\Lambda = M = \{\pm n^\alpha\colon n\in\mathbb N\}$ is a {\sf UP} for the Schwartz space $\mathcal S$
if $\alpha < 1 - \tfrac{\sqrt 2}2$, but their techniques were not strong enough to cover the full range
$\alpha<\tfrac12$.

\subsection{Supercritical and subcritical pairs}

\begin{definition}[{\sf supercritical and subcritical pairs}]\label{def-sub/super-crit}
Let the sequences
$\Lambda=(\lambda_j)_{j\in\Z}$ and $M = (\mu_j)_{j\in\Z}$ be ordered
so that
\begin{align*}
\ldots < \lambda_{j-1} < \lambda_j < \lambda_{j+1} \ldots\,, &\quad
\lim_{j\to\pm\infty} \lambda_j = \pm\infty, \\
\ldots < \mu_{j-1} < \mu_j < \mu_{j+1} \ldots\,, &\quad
\lim_{j\to\pm\infty} \mu_j = \pm\infty\,.
\end{align*}
Given $1<p, q<\infty$ with $\tfrac1p + \tfrac1q =1$, we call a pair
$(\Lambda, M)$ supercritical if
\[
\limsup_{|j|\to \infty} |\lambda_j|^{p-1} \bigl( \lambda_{j+1}-\lambda_j) < \tfrac12,
\ \limsup_{|j|\to \infty} |\mu_j|^{q-1} \bigl( \mu_{j+1}-\mu_j) < \tfrac12\,,
\]
and subcritical if
\[
\liminf_{|j|\to \infty} |\lambda_j|^{p-1} \bigl( \lambda_{j+1}-\lambda_j) > \tfrac12,
\ \liminf_{|j|\to \infty} |\mu_j|^{q-1} \bigl( \mu_{j+1}-\mu_j) > \tfrac12\,.
\]
\end{definition}

\noindent{\bf Example. } Consider the sequences
\[
\Lambda = \bigl\{ \sign(j) \bigl( \tfrac{p}2 a|j| \bigr)^{1/p} \bigr\}_{j\in\Z}, \quad
M = \bigl\{ \sign(j) \bigl( \tfrac{q}2 a|j| \bigr)^{1/q} \bigr\}_{j\in\Z},
\quad \tfrac1p + \tfrac1q = 1, a>0,
\]
and note that
\[
\lim_{|j|\to \infty} |\lambda_j|^{p-1} \bigl( \lambda_{j+1}-\lambda_j) = \tfrac12 a,
\ \lim_{|j|\to \infty} |\mu_j|^{q-1} \bigl( \mu_{j+1}-\mu_j) = \tfrac12 a\,.
\]
Hence, for $a<1$, the pair $(\Lambda, M)$
is supercritical, while for $a>1$, the pair $(\Lambda, M)$ becomes subcritical.
The case $p=q=2$ and $a=1$ corresponds to the pair considered by Radchenko and Viazovska.

\begin{theorem}\label{thm-prelim}
Suppose that $1<p, q<\infty$, $\tfrac1p + \tfrac1q = 1$. Then

\smallskip\noindent {\rm (i)}
any supercritical pair $(\Lambda, M)$ is a uniqueness pair for the Schwartz space $\mathcal S$;

\smallskip\noindent {\rm (ii)}
any subcritical pair $(\Lambda, M)$ is a non-uniqueness pair for the Schwartz space $\mathcal S$.
\end{theorem}
Actually, we prove more. We show that the uniqueness part holds for functions from a
much larger Hilbert space
\[
\mathcal H = \{f\colon f, \widehat{f}\in H_1, \
\| f \|_{\mathcal H}^2 = \| f \|_{H_1}^2 + \| \widehat f \|_{H_1}^2 \}\,,
\]
where $H_1$ is the Sobolev space of functions $f\in L^2(\mathbb R)$ such that
\[
\| f \|_{H_1}^2 = \int_\mathbb R (1+\xi^2)|\widehat{f}(\xi)|^2 \, {\rm d}\xi <\infty\,.
\]
Note that the space $\mathcal H$ is complete and that in this space
the point evaluation is bounded for both $f$ and $\widehat{f}$.
To see that $\mathcal H$ is complete, we let
$(f_n)$ be a Cauchy sequence in $\mathcal H$. Then $(f_n)$ and $(\widehat{f}_n)$
are Cauchy sequences in $H_1$. Hence, $f_n\to f$, $\widehat{f}_n\to g$ in $H_1$,
and therefore, in $L^2(\mathbb R)$. Thus, $\widehat{f}=g$,
and $\widehat{f}_n\to \widehat{f}$ in $H_1$. We conclude that
$f\in\mathcal H$ and $f_n\to f$ in $\mathcal H$.
The boundedness of the point evaluation in $\mathcal H$ is also straightforward:
using the Fourier inversion and the Cauchy--Schwarz inequality, we get
\[
|f(x)| \le
\int_\R \sqrt{1+\xi^2} |\widehat{f}(\xi)|\, \frac{{\rm d}\xi}{\sqrt{1+\xi^2}}
\le \sqrt{\pi} \| f \|_{H_1}\,, \quad x\in\R\,,
\]
and similarly, $|\widehat{f}(\xi)|\le \sqrt{\pi} \| \widehat{f} \|_{H_1}$, $\xi\in\R$.

\medskip
In its turn, the non-uniqueness part of Theorem~\ref{thm-prelim}
requires less stringent assumptions on the pair $(\Lambda, M)$
(see the remark in the beginning of Section~\ref{sect:Non-uniq}) and
holds for a much smaller Gelfand--Shilov space
\[
\mathcal S(p, q) =
\Bigl\{f\colon \exists c=c_f > 0\ {\rm such\ that\ } \int_\mathbb R |f(x)|^2 e^{c|x|^p}\, {\rm d}x < \infty, \
\int_\mathbb R |\widehat{f}(\xi)|^2 e^{c|\xi|^q}\, {\rm d}\xi  < \infty
\Bigr\}\,.
\]
The elements of the Gelfand--Shilov space are entire functions of finite order of growth
($f$ has order $p$ and $\widehat{f}$ has order $q$).

\medskip\noindent{\bf Remark.}
Taking into account the invariance of uniqueness and non-uniqueness pairs under dilations, one can formally
generalize Theorem~\ref{thm-prelim} as follows.

{\em Suppose that, for $1< p, q <\infty$, $\tfrac1p + \tfrac1q = 1$,
\[
\limsup_{|j|\to\infty} |\la_j|^{p-1}(\la_{j+1}-\la_j)=a,
\]
and
\[
\limsup_{|j|\to\infty} |\mu_j|^{q-1}(\mu_{j+1}-\mu_j)=b,
\]
with
\[
a^{1/p}b^{1/q}<\tfrac12\,.
\]
Then $(\Lambda, M)$ is a uniqueness pair.

Similarly, if
\[
\liminf_{|j|\to\infty} |\la_j|^{p-1}(\la_{j+1}-\la_j)=a,
\]
and
\[
\liminf_{|j|\to\infty} |\mu_j|^{q-1}(\mu_{j+1}-\mu_j)=b,
\]
with
\[
a^{1/p}b^{1/q}>\tfrac12\,,
\]
then $(\Lambda, M)$ is a non-uniqueness pair}.

\subsection{Quantitative uniqueness}
To state the rest of our results, we introduce a scale of Hilbert spaces
$\Ho_{s, p, q}$.
\begin{definition}[{\sf Hilbert spaces $\Ho_{s, p, q}$}]
For $p, q > 1$ with $\frac{1}{p} + \frac{1}{q} = 1$ and $s>0$, let $\Ho_{s, p, q}$ be the Hilbert space
\[
\Ho_{s, p, q} = \{f\colon
f\in H_{qs}, \widehat{f}\in H_{ps}, \ \| f \|_{\mathcal H_{s, p, q}}^2 = \| f \|_{H_{ps}}^2 + \| \widehat f \|_{H_{qs}}^2\,
\}\,,
\]
where
$H_t$ is the Sobolev space of functions $f\in L^2(\mathbb R)$ such that
\[
\| f \|_{H_t}^2 = \int_\mathbb R (1+|\xi|^{2t})|\widehat{f}(\xi)|^2 \, {\rm d}\xi <\infty\,.
\]
\end{definition}
\noindent
Note that $\mathcal H_{\frac12, 2, 2} = \Ho$.

We will use the spaces $\mathcal H_{s, p, q}$ only for the values of $s$ such that $s\min(p, q) \ge 1$ so that $\mathcal H_{s, p, q}\subset \mathcal{H}$.

We also need to introduce the notion of $p$-separated sequences, that is the sequences for which the consecutive elements can not be too close.
\begin{definition}[{\sf $p$-separated sequences}]\label{def:ser_seq}
Let $\Lambda=(\lambda_j)_{j\in\Z}$ be a sequence ordered
so that
$$
\ldots < \lambda_{j-1} < \lambda_j < \lambda_{j+1} \ldots\,, \quad
\lim_{j\to\pm\infty} \lambda_j = \pm\infty. \\
$$
We call $\Lambda$ a $p$-separated sequence if there exists $c > 0$ such that for all $j\in \Z$ we have
$$\lambda_{j+1}-\lambda_j \ge \frac{c}{(1+\min(|\lambda_j|, |\lambda_{j+1}|))^{p-1}}\,.$$
\end{definition}

We will tacitly use the following simple claim:
\begin{claim}\label{claim:separation}
For any supercritical pair $(\Lambda, M)$,
there exists a supercritical pair $(\Lambda', M')$ such that $\Lambda' \subset \Lambda$, $M'\subset M$, and $\Lambda'$ is $p$-separated while $M'$ is $q$-separated.
\end{claim}

\noindent{\em Proof:}
We only show how to select the needed subsequence $\Lambda'_+$ from
$ \Lambda \cap (0, +\infty) = \{\lambda_k, \lambda_{k+1}, \ldots\}$, the other three cases are handled similarly. Let
$$a = \limsup |\lambda_j|^{p-1}(\lambda_{j+1}-\lambda_j) < \frac{1}{2}$$
and consider a small number $\delta > 0$. We build the set $\Lambda'_+$ iteratively, starting with $\Lambda'_+ = \{ \lambda_k\}$. We will be looking at the elements of $\Lambda\cap (0, +\infty)$ in order. Assume that we are currently looking at $\lambda_m$ and let $\lambda_n$ be the last element that we added to $\Lambda'_+$. We add $\lambda_m$ to $\Lambda'_+$ if and only if $\lambda_n^{p-1}(\lambda_m - \lambda_n) \ge \delta$. Clearly, the resulting set $\Lambda'_+$ is $p$-separated.
We show that if $a+2\delta < \frac12$, then this set remains $p$-supercritical.

Let us enumerate $\Lambda'_+ = \{\gamma_1, \gamma_2, \ldots\}$ and show that, for sufficiently
large $l$, we have $\gamma_l^{p-1}(\gamma_{l+1}-\gamma_l) < a+2\delta$.
Assume that in the original sequence $\gamma_l = \lambda_n$, $\gamma_{l+1} = \lambda_m$.
If $m = n+1$ then $\gamma_l^{p-1} (\gamma_{l+1}-\gamma_l) < a+\delta$, provided that
$l$ is large enough.

Now, assume that $m>n+1$, that is, $\lambda_{m-1}\notin \Lambda'_+$. Then
\begin{multline*}
\gamma_l^{p-1}(\gamma_{l+1}-\gamma_l) = \la_n^{p-1}(\lambda_{m}-\lambda_{m-1})
+ \la_n^{p-1}(\lambda_{m-1}-\lambda_n) \\
< \la_{m-1}^{p-1}(\lambda_{m}-\lambda_{m-1})
+ \la_n^{p-1}(\lambda_{m-1}-\lambda_n)
< (a+\delta) + \delta = a + 2\delta\,,
\end{multline*}
provided again that $l$, and therefore $m$, are large enough. \hfill $\Box$

\medskip
Thus, there is no loss in generality to assume in Theorem \ref{thm-prelim} that $\Lambda$ and $M$ are $p$- and $q$-separated, respectively. We will call such pairs separated with parameters $p$ and $q$.

\medskip
Our next result is a quantitative version of the uniqueness part of Theorem~\ref{thm-prelim}.
\begin{theorem}\label{thm-stable-uniq}
Let $(\Lambda, M)$ be a separated supercritical pair with parameters
$p, q>1$, $\tfrac1p + \tfrac1q = 1$,
and let $s>0$ be such that $s\min(p, q) \ge 1$. Then there exist numbers $c, C > 0$, such that,
for every $f\in\mathcal H$,
\begin{multline*}
c\| f \|_{\mathcal H_{s, p, q}}^2 \le \,
\Bigl[\,
\sum_{\lambda\in\Lambda} (1+|\lambda|)^{(2s-1)p + 1} |f(\lambda)|^2 \\
+ \sum_{\mu\in M} (1+|\mu|)^{(2s-1)q+1} |\widehat{f}(\mu)|^2
\,\Bigr]\,\le C \| f \|_{\mathcal H_{s, p, q}}^2.
\end{multline*}
\end{theorem}
Note that the sums in the middle are comparable to the Riemann sums for
the integrals that define $||f||^2_{\Ho_{s, p, q}}$. By the left-hand inequality,
if the expression in the brackets is finite for a function $f\in \Ho$, then $f$ belongs to the space $\Ho_{s, p, q}$.

\section{Corollaries and related work}

\subsection{Fast decay along a supercritical pair yields membership in $\mathcal S$}
Combining Theorem~\ref{thm-stable-uniq} and Claim~\ref{claim:separation} with a simple fact that
$\bigcap_{s>0} \Ho_{s, p, q} = \mathcal{S}$, we get the following corollary:
\begin{corollary}\label{schwa theorem}
Let $(\Lambda, M)$ be a supercritical pair. The function $f\in \Ho$ belongs to the Schwartz space $\mathcal{S}$ if and only if for all $r<\infty$ we have
$$|f(\lambda_j)| = O\left( |\lambda_j|^{-r} \right), \quad |\widehat{f}(\mu_j)|
=  O\left( |\mu_j|^{-r} \right),
\quad |j|\to\infty\,.
$$
\end{corollary}

\subsection{Interpolation formula }
Theorem \ref{thm-stable-uniq} essentially says that suitably normalized reproducing kernels at the sets $\Lambda$ and $M$ form a frame for the Hilbert space $\mathcal{H}_{s, p, q}$. It is well-known that any such frame bound gives rise to an interpolation formula, see e.g. \cite[Section~5.1]{Grochenig} or \cite[Section~1.8]{OUbook}. For the reader's convenience and to keep this text self-contained,
we will outline the argument in our case following~\cite{OUbook}.

\begin{corollary}\label{thm-interpol}
Let $(\Lambda, M)$ be a supercritical pair with parameters $p, q>1$, $\tfrac1p + \tfrac1q = 1$,
and let $s>0$ satisfy $s\min(p, q)\ge 1$. Then there exist functions $(a_\lambda)_{\lambda\in\Lambda}$ and $(b_\mu)_{\mu\in M}$ in $\mathcal H_{s, p, q}$, such that, for every $f\in\mathcal H_{s, p, q}$, we have
\begin{equation}\label{eq-interpol}
f = \sum_{\lambda\in\Lambda} f(\lambda) a_\lambda
+ \sum_{\mu\in M} \widehat{f}(\mu) b_\mu\,
\end{equation}
and the series on the RHS converges in $\mathcal H_{s, p, q}$.

Furthermore, the functions $a_\lambda$ and $b_\mu$ can be chosen so that
\[
\| a_\la \|_{\mathcal H_{s, p, q}} \le C(1+|\la|)^{(s-\frac12)p+\frac12}, \quad \la\in\Lambda,
\]
and
\[
\| b_\mu \|_{\mathcal H_{s, p, q}} \le C(1+|\mu|)^{(s-\frac12)q+\frac12}, \quad \mu\in M.
\]
\end{corollary}
Since $|f(x)| \le \sqrt{\pi}\| f \|_{H_1}$ and
$\| f \|_{H_1} \le \| f \|_{\Ho} \le \sqrt{2}\| f \|_{\Ho_{s,p,q}}$,
we have $|f(x)| \le \sqrt{2\pi}\| f \|_{\Ho_{s,p,q}}$,
uniformly in $x\in\mathbb R$, and similarly
$|\widehat{f}(\xi)| \le \sqrt{2\pi} \| f \|_{\Ho_{s,p,q}}$ uniformly in $\xi\in\mathbb R$,
so the series on the RHS of~\eqref{eq-interpol} converges  uniformly
both on the space and the Fourier sides.
Furthermore, by a similar argument, if $s\min (p, q)\ge n+1$, $n\in \N$,
then the series in~\eqref{eq-interpol} converges uniformly together with the first $n$ derivatives
both on the space and the Fourier sides.

On the other hand, it is important to note that the functions $a_\lambda, b_\mu$ depend on $s$ and that
we do not know whether there exists a similar interpolation formula convergent in the topology of the Schwartz space $\mathcal S$.

\medskip\noindent{\em Proof of Corollary~\ref{thm-interpol}}:
By Claim~\ref{claim:separation}, we can assume that the pair $(\Lam, M)$ is separated.

Let $\ell^2(\Lambda, M)$ be the Hilbert space of (complex-valued) sequences $x$ indexed by
$\Lambda \cup M$ such that
\[
\| x \|^2\, \stackrel{{\rm def}}=\, \sum_{\la\in\Lambda} |x(\la)|^2 (1+|\la|)^{(2s-1)p +1}
+ \sum_{\mu\in M} |x(\mu)|^2 (1+|\mu|)^{(2s-1)q +1} < \infty\,.
\]
Then the map
\[ f\mapsto \bigl(\, (f(\la))_{\la\in\Lambda}, (\widehat{f}(\mu))_{\mu\in M}
\, \bigr) \]
defines a linear operator \[ T\colon \mathcal H_{s, p, q} \to \ell^2(\Lambda, M), \]
and, by Theorem~\ref{thm-stable-uniq}, there exists a positive $C$ such that, for each $f\in
\mathcal H_{s, p, q}$,
\begin{equation}\label{eq:quasi-isom}
C^{-1} \| f \|_{\mathcal H_{s, p, q}} \le \| Tf \|_{\ell^2(\Lambda, M)}
\le C \|f\|_{\mathcal H_{s, p, q}}\,.
\end{equation}
Let $X=T\mathcal H_{s, p, q}$ be the image of the space $\mathcal H_{s, p, q}$.
By~\eqref{eq:quasi-isom}, this is a closed subspace of $\ell^2(\Lambda, M)$. Denote by
$P\colon \ell^2(\Lambda, M)\to X$ the orthogonal projection to $X$.

Let $\delta_\la$, $\la\in\Lambda$, and $\delta_\mu$, $\mu\in M$, be the ``delta-functions''
in $\ell^2(\Lambda, M)$, i.e., $\delta_\la(\lambda)=1$, $\delta_\la(\lambda')=\delta_\la(\mu)=0$
for $\lambda'\in\Lambda\setminus\{\la\}$, $\mu\in M$, and similarly for $\delta_{\mu}$.
Let $e_\la = P\delta_\la$, $e_\mu = P\delta_\mu$. Then, for each $x=Tf$ in $X$, we have
\[
x = \sum_{\la\in\Lambda} x(\la) e_\la + \sum_{\mu\in M} x(\mu) e_\mu
= \sum_{\la\in\Lambda} f(\la) e_\la + \sum_{\mu\in M} \widehat{f}(\mu) e_\mu\,,
\]
with the series convergent in $\ell^2(\Lambda, M)$. Applying the inverse operator
$T^{-1}$ (which is, by~\eqref{eq:quasi-isom}, a bounded operator from $X$ to $\mathcal H_{s, p, q}$), letting $a_\la = T^{-1}e_\la$, $b_\mu=T^{-1}e_\mu$, and noting that
\[
\| a_\la \|_{\mathcal H_{s, p, q}} \le \| T^{-1} \|_{X\to \mathcal H_{s, p, q} }\,
\|e_\la\|_{\ell^2(\Lambda, M)} \le C \| \delta_\la \|_{\ell^2(\Lambda, M)}
= C(1+|\la|)^{(s-\frac12)p+\frac12}, \quad \la\in\Lambda,
\]
and, similarly, $\| b_\mu \|_{\mathcal H_{s, p, q}} \le C(1+|\mu|)^{(s-\frac12)q+\frac12}$,
$\mu\in M$, we conclude the proof. \hfill $\Box$

\subsection{Crystalline measures}
\medskip Keeping $x\in\R\setminus\Lambda$ fixed, we can rewrite
the interpolation formula~\eqref{eq-interpol} in the form
\[
f(x) - \sum_{\lambda\in\Lambda} f(\lambda) a_\lambda (x)
= \sum_{\mu\in M} \widehat{f}(\mu) b_\mu (x),
\qquad f\in\mathcal S.
\]
Then (keeping in mind that $\sup_\R |a_\lambda| \le C \| a_\lambda \|_{\Ho_{s, p, q}}$ and
$\sup_\R |b_\mu| \le C \| b_\mu \|_{\Ho_{s, p, q}}$),
we see that the locally finite measure $\nu_x = \sum_{\mu\in M} b_\mu(x)\delta_\mu$
is a tempered distribution supported on the set $M$, and its distributional
Fourier transform equals to $\widehat{\nu}_x =  \delta_x - \sum_{\lambda\in\Lambda} a_\lambda(x)\delta_\lambda  $, which
is a locally finite measure supported on the set $\Lambda\bigcup \{ x\}$.
Applying this argument with some $\lambda'\in\Lambda$, $\Lambda' = \Lambda\setminus\{\lambda'\}$
in place of $\Lambda$, and $x=\lambda'$, we arrive at a curious conclusion:
\begin{corollary}\label{cor-crystalline}
For any supercritical pair $(\Lambda, M)$, there exists a non-trivial locally finite measure $\nu$
supported by $M$ that is a tempered distribution, whose distributional Fourier transform $\widehat\nu$ is a locally finite measure supported by $\Lambda$.
\end{corollary}
Such measures $\nu$ are called {\em one-dimensional quasicrystals} or {\em crystalline measures} according to Dyson~\cite{Dyson} and Meyer~\cite{Meyer}.

\medskip
In fact, we can prove slightly more. For a pair of discrete sets $\Lambda$, $M$, we denote by $\mathcal{N}(\Lambda, M)$ the vector space of locally finite measures $\nu$ supported on $M$ which are tempered distributions, whose distributional Fourier transforms $\widehat\nu$ are locally finite measures supported on $\Lambda$. Then, since removing finitely many points from a supercritical pair leaves a supercritical pair, the space $\mathcal{N}(\Lambda, M)$ is not only non-trivial,
but also infinite-dimensional for any supercritical pair $(\Lambda, M)$. On the other hand, for the subcritical pair, this space is always finite-dimensional, see Section~\ref{subsect-FreeInterpol}.

\subsection{Related work}
There are several recent results inspired by the aforementioned
Radchenko--Viazovska paper~\cite{Rad-Via}.
Bondarenko, Radchenko and Seip~\cite{BRS} revealed a collection of similar interpolation formulas,
the most curious one corresponding to $\Lambda = \{ (4\pi)^{-1} \log n\colon n\in\mathbb N\}$ and
$M = \{ i (\rho - \tfrac12) \}$, where $\rho$ runs through all non-trivial zeroes of the
Riemann zeta function with positive imaginary part (assuming the Riemann hypothesis and simplicity of
zeroes). Similarly to~\cite{Rad-Via}, they used an arithmetic approach based on the theory of modular forms and Dirichlet series. In~\cite{Ram-Sousa2}, Ramos and Sousa
managed to slightly perturb the nodes of the interpolation formula found in~\cite{Rad-Via}.

Recently, in~\cite{AA} Adve found conditions for the pair of discrete set $(\Lambda, M)\subset R^d \times R^d$ to be a UP/NUP for the Schwartz space $S(R^d)$. His conditions are less precise than the ones supplied by Theorem~\ref{thm-prelim} in the case $d=1$ (apart from that, they are of uniform nature, like in our toy-model Theorem~1A, rather than asymptotic). We note that  his construction of NUPs does not appeal to complex analysis.  It should also be mentioned that other, denser than in~\cite{AA}, multi-dimensional discrete UP were recently constructed by Ramos and Stoller~\cite{R-S} and by Viazovska~\cite{V}.

In~\cite{Kulikov}, the
first-named author proved that interpolation formulas like~\eqref{eq-interpol} can only exist if the following inequality for the counting functions of the sets $\Lambda, M$ is satisfied for all $W, T > 1$:
\begin{equation*}
|\Lambda \cap [-T, T]| + |M \cap [-W, W]| \ge 4WT - C\log^2(4WT)\,.
\end{equation*}
One can easily see that the supercriticality condition in Theorem~\ref{thm-interpol} implies this inequality (even with $4$ replaced by some bigger constant), while, for subcritical pairs
in the non-uniqueness part of Theorem~\ref{thm-prelim}, there exist values of $T, W$ that violate it regardless of the value of the constant $C$.

For the reader interested in one-dimensional quasicrystals, we mention recent works by Kurasov and Sarnak~\cite{KS}, Lev and Reti~\cite{LP1}, Olevskii and Ulanovskii~\cite{OU}, and Meyer~\cite{Meyer}, which also contain further references.

\section{Poincar\'e-Wirtinger-type inequalities}\label{sect:Wirt}
The central role in our arguments will be played by the Poincar\'e--Wirtinger inequality
and its extensions. The most classical version states that if $f\in C^1[0, 1]$
and $f(0)=f(1)=0$, then\footnote{
For the proof, extend the function $f$ to the odd $2$-periodic function $\widetilde f$
and let $g(t)=\widetilde{f}(2t)$. The function $g$ is a $1$-periodic $C^1$-function
with zero mean. Expanding $g$ and $g'$ into the Fourier series and using that
$\widehat{g'}(k)= 2\pi {\rm i} k\widehat{g}(k)$, $k\in\mathbb Z$, we get
\[
4 \int_0^1 |f'|^2 = \int_{\mathbb R/\mathbb Z} |g'|^2
= 4\pi^2 \sum_{k\ne 0} k^2 |\widehat{g}(k)|^2
\ge 4\pi^2 \sum_{k\ne 0} |\widehat{g}(k)|^2
= 4\pi^2 \int_{\mathbb R/\mathbb Z} |g|^2
= 4\pi^2 \int_0^1 |f|^2\,.
\]
}
\begin{equation*}\label{eq:PW}
\int_0^1 |f|^2 \le \Bigl( \frac1{\pi} \Bigr)^2\, \int_0^1 |f'|^2\,.
\end{equation*}

We start with a ``stable version'' of this inequality, which does not assume that $f$ vanishes
at the endpoints.

\begin{lemma}\label{wirt1}
For any $\eps > 0$, all $a < b$ and all functions
$f\in H_1(\R)$, we have
\[
\int_a^b |f|^2 \\
\le (1+\eps)\, \left(\frac{b-a}{\pi}\right)^2\,
\int_a^b |f'|^2\, + (1+\eps^{-1}) (b-a)(|f(a)|^2 + |f(b)|^2).
\]
\end{lemma}

\noindent{\em Proof:}
Since  smooth functions are dense in $H_1(\R)$ and the point evaluation is bounded in
$H_1(\R)$, it suffices to prove the estimate for a smooth function $f$ on $[a, b]$.

By the scale- and shift-invariance we can assume that $a = 0$ and $b = 1$.
Consider the function $g = f - \ell $, where the linear function $\ell$ is taken so that
$\ell(0)=f(0)$ and $\ell (1) = f(1)$, i.e., $g(0)=g(1)=0$. Then, by the classical Poincar{\'e}--Wirtinger inequality,
\[
\int_0^1 |g|^2 \le \frac{1}{\pi^2}\, \int_0^1 |g'|^2\,.
\]
Furthermore, by the Cauchy--Scwarz inequality,
\[
\int_0^1 |f|^2 = \int_0^1 |g+\ell|^2 \le
(1+\eps)\int_0^1 |g|^2 + (1+\eps^{-1})\int_0^1|\ell|^2\,,
\]
and
\[
\int_0^1 |f'|^2 = \int_0^1 |g'|^2 + \int_0^1 |\ell'|^2 \ge \int_0^1 |g'|^2
\]
(since the function $g'$ has zero mean and, therefore, is orthogonal to constant functions).
At last, $|\ell| \le \max (|\ell (0)|, |\ell (1)|) = \max (|f(0)|, |f(1)|)$ everywhere on $[0, 1]$, whence
\[
\int_0^1 |\ell|^2 \le \max (|f(0)|^2, |f(1)|^2) \le |f(0)|^2 + |f(1)|^2\,.
\]
Combining these estimates, we complete the proof. \hfill $\Box$

We also need the following trace-type bound:
\begin{lemma}\label{trace-claim}
For all $a<b$ and all $f\in H_1(\mathbb R)$, we have
$$
\frac{1}{b-a} |f(a)|^2
\le \frac2{(b-a)^2}\int_a^b |f|^2 + \frac23 \int_a^b |f'|^2\,.
$$
\end{lemma}

\noindent{\em Proof:}
As above, we assume that $a = 0$, $b = 1$, and that the function $f$ is smooth.
Integrating the identity
\[
f(0) = f(x) - \int_0^x f'
\]
over $x\in [0, 1]$, we get
\[
f(0) = \int_0^1 f - \int_0^1 (1-x)f'\,,
\]
whence,
\[
|f(0)| \le \int_0^1 |f| + \int_0^1 (1-x)|f'|\,,
\]
and
\begin{multline*}
|f(0)|^2 \le 2\Bigl[ \Bigl( \int_0^1 |f| \Bigr)^2 + \Bigl( \int_0^1 (1-x)|f'| \Bigr)^2 \Bigr] \\
\le 2\int_0^1 |f|^2 + 2 \int_0^1 (1-x)^2 \, \cdot \int_0^1 |f'|^2
= 2\int_0^1 |f|^2 + \frac23\, \int_0^1 |f'|^2\,,
\end{multline*}
proving the lemma. \hfill $\Box$

\begin{definition}[{\sf $l$-dense sets}]\label{def:dense}
A discrete set $\Gamma\subset\mathbb R$ is $l$-dense if $\mathbb R\setminus\Gamma$ contains
no interval of length greater than $l$.
\end{definition}

\begin{lemma}\label{wirt2} \mbox{}
\smallskip\noindent {\rm (i)}
Let
$t>0$, $0<\eps<1$, and let $\Gamma$ be a $(1-\eps)(2t)^{-1}$-dense discrete subset
of $\mathbb R$.
There exists a constant $C_\eps$ such that, for all convex increasing $C^1$-functions
$\Phi\colon [0, \infty) \to \mathbb R$ and all functions $f\in\Ho$, we have
\[
\Phi(t^2)\, \int_{\mathbb R} |f(x)|^2\, {\rm d}x
\le \int_{\mathbb R} \Phi(\xi^2) |\widehat{f}(\xi)|^2\, {\rm d}\xi
+ C_\eps t\Phi'(t^2)\, \sum_{\gamma\in\Gamma} |f(\gamma)|^2\,.
\]

\noindent {\rm (ii)} Let $t>0$, and let $\Gamma$ be a $(2t)^{-1}$-dense discrete subset
of $\mathbb R$.
Then, for all convex increasing $C^1$-functions
$\Phi\colon [0, \infty) \to \mathbb R$ and all functions $f\in\Ho$ vanishing
on $\Gamma$, we have
\[
\Phi(t^2)\, \int_{\mathbb R} |f(x)|^2\, {\rm d}x
\le \int_{\mathbb R} \Phi(\xi^2) |\widehat{f}(\xi)|^2\, {\rm d}\xi\,.
\]
\end{lemma}
The case $\Phi(t)=t^\theta$, $\theta\ge 1$, will be of a particular interest for us later.

\medskip\noindent{\em Proof:}
First, we consider the case $\Phi(t) = t$. Summing the estimates from Lemma~\ref{wirt1} over all gaps between the  points of the set $\Gamma$ and using that the length of the gaps cannot exceed
$(1-\eps)/(2t)$, we get
\[
\int_\mathbb R |f|^2 \le
\frac{(1+\eps)(1-\eps)^2}{(2\pi t)^2}\, \int_\mathbb R |f'|^2\,
+ \frac{(1+\eps^{-1})(1-\eps)}{t}\, \sum_{\gamma\in\Gamma} |f(\gamma)|^2\,,
\]
whence (using the Plancherel theorem),
\[
t^2 \int_{\mathbb R} |f|^2
\le  \int_\mathbb R \xi^2 |\widehat{f}(\xi)|^2\, {\rm d}\xi
+ C_\eps t \sum_{\gamma\in\Gamma} |f(\gamma)|^2\,.
\]

For the general case, we rewrite the case $\Phi(t) = t$ as
$$\int_{|\xi|\le t} (t^2-\xi^2) |\widehat{f}(\xi)|^2\, {\rm d}\xi
\le \int_{|\xi| > t} (\xi^2-t^2) |\widehat{f}(\xi)|^2\, {\rm d}\xi
+ C_\eps t \sum_{\gamma\in \Gamma} |f(\gamma)|^2.$$
Multiplying this inequality by $\Phi'(t^2)$, we get
\begin{multline*}
\Phi'(t^2) \int_{|\xi|\le t} (t^2-\xi^2) |\widehat{f}(\xi)|^2\, {\rm d}\xi \\
\le \Phi'(t^2) \int_{|\xi| > t} (\xi^2-t^2) |\widehat{f}(\xi)|^2\, {\rm d}\xi
+ C_\eps t \Phi'(t^2) \sum_{\gamma\in \Gamma} |f(\gamma)|^2.
\end{multline*}
Noting that, for $|\xi| > t$, we have $\Phi(\xi^2)-\Phi(t^2)\ge \Phi'(t^2)(\xi^2-t^2)$, and that,
for $|\xi| \le t$, we have $\Phi(t^2)-\Phi(\xi^2) \le \Phi'(t^2)(t^2-\xi^2)$, we get
\begin{multline*}
\int_{|\xi|\le t} |\widehat{f}(\xi)|^2 (\Phi(t^2)-\Phi(\xi^2))\, {\rm d}\xi  \\
\le \int_{|\xi| > t} |\widehat{f}(\xi)|^2 (\Phi(\xi^2)-\Phi(t^2))\, {\rm d}\xi
+ C_\eps t \Phi'(t^2)\sum_{\gamma\in \Gamma} |f(\gamma)|^2\,,
\end{multline*}
which after rearranging gives us part (i) of the lemma.

To prove the second part, we notice that the set $\Gamma$ is $(1-\eps)(2t')^{-1}$-dense
with any $t'<t$ and $\eps<1-t^{-1}t'$. Applying part (i), recalling that $f$ vanishes on $\Gamma$,
and letting $t'\uparrow t$, we get part (ii). \hfill $\Box$

\section{Toy model: the Poincar\'e--Wirtinger inequality in action}
Here we prove the uniqueness part of Theorem~\ref{thm-prelim}, as well as
Theorem~\ref{thm-stable-uniq}, in the case $p = q = 2$, $s=\tfrac12$,
under a more restrictive {\em uniform supercriticality assumption}. In this case,
the proofs of these results become short and free of cumbersome technical details.

\subsection{Uniformly supercritical uniqueness pairs}
As before, we order the sequences $\Lambda=(\lambda_j)_{j\in\Z}$ and $M = (\mu_j)_{j\in\Z}$,
so that
\begin{align*}
\ldots < \lambda_{j-1} < \lambda_j < \lambda_{j+1} \ldots\,, &\quad
\lim_{j\to\pm\infty} \lambda_j = \pm\infty, \\
\ldots < \mu_{j-1} < \mu_j < \mu_{j+1} \ldots\,, &\quad
\lim_{j\to\pm\infty} \mu_j = \pm\infty\,.
\end{align*}
\begin{definition}[{\sf uniform supercriticality}]\label{def-uniform-supercrit}
We call a pair $(\Lambda, M)$ uniformly supercritical, if
\begin{equation}\label{eq-uniform-supercit}
\sup_{j\in\Z}\, \max\{ |\lambda_j|, |\lambda_{j+1}| \} \cdot
 (\lambda_{j+1}-\lambda_j), \
\sup_{j\in\Z}\, \max\{ |\mu_j|, |\mu_{j+1}| \} \cdot
 (\mu_{j+1}-\mu_j) <  \tfrac12\,.
\end{equation}
\end{definition}

\medskip\noindent{\bf Theorem 1A.}
{\em Let $(\Lambda, M)$ be a uniformly supercritical pair.
Then the only function $f\in \Ho$ such that $f(\lambda_j) = \widehat{f}(\mu_j) = 0$, $j\in \Z$,
is $f= 0$}.

\medskip\noindent{\em Proof:}
Choose $a<1$ so that the LHS of~\eqref{eq-uniform-supercit}
is less than $\tfrac12\, a$. Then, we have
\begin{align*}
\int_\R x^2|f(x)|^2\, {\rm d}x &= \sum_{j\in \Z}\,
\int_{\lambda_j}^{\lambda_{j+1}}x^2|f(x)|^2\, {\rm d}x
\\
&\le
\sum_{j\in \Z}\,  \max (|\lambda_j|, |\lambda_{j+1}|)^2\,
\int_{\lambda_j}^{\lambda_{j+1}} |f(x)|^2\, {\rm d}x
\\
&\le
\sum_{j\in \Z}\, \max (|\lambda_j|, |\lambda_{j+1}|)^2\,
\left(\frac{\lambda_{j+1}-\lambda_j}{\pi}\, \right)^2\,
\int_{\lambda_j}^{\lambda_{j+1}} |f'(x)|^2\, {\rm d} x
\\
&\le
\frac{a^2}{(2\pi)^2}\, \sum_{j\in \Z}\, \int_{\lambda_j}^{\lambda_{j+1}} |f'(x)|^2\, {\rm d}x
\\
&= \frac{a^2}{(2\pi)^2}\, \int_\R |f'(x)|^2 \, {\rm d} x
\\
&= a^2 \int_\R y^2 |\widehat{f}(y)|^2\, {\rm d}y,
\end{align*}
where in the second inequality we used the classical
Poincar{\'e}--Wirtinger inequality.

Similarly, we get
$$\int_\R y^2|\widehat{f}(y)|^2\, {\rm d}y \le a^2 \int_\R x^2|f(x)|^2\, {\rm d}x\,. $$
Since $a <1 $, this implies that $\displaystyle \int_\R x^2|f(x)|^2\, {\rm d}x = 0$, i.e.,
$f= 0$ almost everywhere. \hfill $\Box$

\medskip
By being slightly more careful, we can assume in Theorem~1A that the LHS of~\eqref{eq-uniform-supercit} are only at most $\tfrac{1}{2}$ and not strictly smaller than $\frac{1}{2}$. As the Radchenko--Viazovska construction~\cite{Rad-Via} shows, this is no longer true in the case
of more general Theorem~\ref{thm-prelim}.

Note that it is enough to add finitely many points to each of the sets $\Lambda$ and $M$,
to turn any supercritical pair $(\Lambda, M)$ (with $p=q=2$)
into the uniformly supercritical one.
Therefore, if $(\Lambda, M)$ is a supercritical pair with $p=q=2$, the linear space
$\{ f\in\mathcal H\colon f|_\Lambda = 0, \widehat{f}|_M = 0\}$
can be only finite-dimensional.
Unfortunately, we did not succeed in finding a short way to conclude from here
that this linear space contains only the identically vanishing function.
Nevertheless, the same idea based on the Poincar\'e--Wirtinger inequality
will play a crucial role in the proof of the uniqueness part of Theorem~\ref{thm-prelim}.

\subsection{Quantitative uniqueness}
The next result is a uniform version of Theorem~\ref{thm-stable-uniq}.

\medskip\noindent{\bf Theorem~2A.}
{\em
Let $(\Lambda, M)$ be a separated uniformly supercritical pair with $p=q=2$.
There exist $A, B >0$ such that, for all $f\in \Ho$, we have}
\begin{equation}\label{frame unif}
A\, ||f||_\Ho^2 \le \sum_{\lambda\in\Lambda} (1+|\lambda|)|f(\lambda)|^2 +
\sum_{\mu\in M}(1+|\mu|)|\widehat{f}(\mu)|^2\le B\, ||f||_\Ho^2.
\end{equation}
We start with the left-hand part of~\eqref{frame unif}, which
does not need the separation assumption.

\medskip\noindent{\em Proof of the left-hand part of~\eqref{frame unif}:}
Choose $a<1$ so that the LHS of~\eqref{eq-uniform-supercit}
is less than $\tfrac12\, a$.
Repeating the proof of Theorem~1A, with Lemma~\ref{wirt1} in place of
the Poincar{\'e}--Wirtinger inequality, we get
\begin{multline}\label{eq:new1}
\int_\R x^2|f(x)|^2\, {\rm d}x \le a^2(1+\eps)\,\int_\R \xi^2|\widehat{f}(\xi)|^2\, {\rm d}\xi
\\
+ 2(1+\eps^{-1}) \sum_{j\in\mathbb Z} \max(|\la_{j-1}|, |\la_j|, |\la_{j+1}|)^2 (\la_{j+1}-\la_j) |f(\la_j)|^2\,.
\end{multline}
By the uniform supercriticality,
$\max(|\la_{j-1}|, |\la_j|, |\la_{j+1}|) \le C (1+|\la_j|)$, and therefore, the second term on the RHS of~\eqref{eq:new1} is bounded by
\[
C_\eps\, \sum_{j\in\mathbb Z}\, \max(|\la_{j-1}|, |\la_j|, |\la_{j+1}|) |f(\lambda_j)|^2
\le C_\eps\, \sum_{\la\in\Lambda} (1+|\la|) |f(\la)|^2\,.
\]
Hence,
\[
\int_\R x^2|f(x)|^2\, {\rm d}x \le a^2(1+\eps)\,\int_\R \xi^2|\widehat{f}(\xi)|^2\, {\rm d}\xi
+ C_\eps\, \sum_{\la\in\Lambda} (1+|\la|) |f(\la)|^2\,.
\]
Similarly,
\begin{equation*}
\int_\R \xi^2|\widehat{f}(\xi)|^2\, {\rm d}\xi
\le a^2(1+\eps)\, \int_\R x^2|f(x)|^2\, {\rm d}x + C_\eps\sum_{\mu\in M} (1+|\mu|) |\widehat{f}(\mu)|^2.
\end{equation*}
Adding these estimates, and choosing $\eps > 0$ so that $a^2(1+\eps) < 1 - \delta$, $\delta > 0$,
we get
\begin{multline*}
\delta \Bigl[ \int_\R x^2|f(x)|^2\, {\rm d}x + \int_\R \xi^2|\widehat{f}(\xi)|^2\, {\rm d}\xi \Bigr]
\\
\le C_\eps\, \Bigl[
\sum_{\lambda\in\Lambda}\, (1+|\lambda|) |f(\lambda)|^2 +
\sum_{\mu\in M} \, (1+|\mu|)|\widehat{f}(\mu)|^2 \Bigr].
\end{multline*}
Recalling that
$ \| f \|_{\Ho}^2 $ is dominated by $ \displaystyle \int_{\mathbb R} x^2 |f|^2 +
\int_{\mathbb R} \xi^2 |\widehat{f}|^2 $
(for instance, by the classical Heisenberg inequality\footnote{
We use it in the form
\[
(2\pi)^{-1} \| f \|^2_{L^2(\mathbb R)} \le
\int_{\mathbb R} x^2 |f(x)|^2\, {\rm d}x + \int_{\mathbb R} \xi^2
|\widehat{f} (\xi)|^2\, {\rm d}\xi\,.
\]
}~\cite[Section~2.9]{DM}),
we get
\begin{equation*}
\delta \| f \|_{\Ho}^2 \le \widetilde{C}_\eps\, \Bigl[
\sum_{\lambda\in\Lambda}\, (1+|\lambda|) |f(\lambda)|^2 +
\sum_{\mu\in M} \, (1+|\mu|)|\widehat{f}(\mu)|^2 \Bigr],
\end{equation*}
which implies the left-hand part of estimate~\eqref{frame unif} with
$A = \delta\, \widetilde{C}_\eps^{-1}$. \hfill $\Box$

\medskip\noindent{\em Proof of the right-hand part of~\eqref{frame unif}:}
Here we will use only the separation condition, choosing
$c$ so that $\lambda_{j+1}-\lambda_j \ge c (1+\min(|\lambda_j|, |\lambda_{j+1}|))^{-1}$,
$j\in\mathbb Z$. For $\la_j\ge 0$, applying Lemma~\ref{trace-claim} with
$a=\lambda_j$, $b=\lambda_j + c(1+\lambda_j)^{-1}$ so that $[a, b] \subset [\la_j, \la_{j+1}]$,
we get
\begin{align*}
c^{-1}(1+\la_j) |f(\la_j)|^2 &\le 2 c^{-2}(1+\la_j)^2\,
\int_{\la_j}^{\la_j +c (1+\la_{j})^{-1}} |f|^2
+ \frac23\, \int_{\la_j}^{\la_j +c (1+\la_{j})^{-1}} |f'|^2 \\
&\le 2 c^{-2}(1+\la_j)^2\, \int_{\la_j}^{\la_{j+1}} |f|^2 + \frac23\,  \int_{\la_j}^{\la_{j+1}} |f'|^2 \\
&\le 2 c^{-2}\, \int_{\la_j}^{\la_{j+1}}(1+x^2) |f|^2 + \frac23\, \int_{\la_j}^{\la_{j+1}} |f'|^2\,.
\end{align*}
For $\la_j<0$, we apply Lemma~\ref{trace-claim} to the function $f(-x)$
with $a=-\lambda_j$, $b=-\lambda_j + c(1+|\lambda_j|)^{-1}$ (note that now $[a, b]\subset [-\la_j, -\la_{j-1}]$) and get
\[
c^{-1}(1+|\la_j|) |f(\la_j)|^2
\le 2 c^{-2}\,\int_{\la_{j-1}}^{\la_j}(1+x^2) |f|^2 + \frac23\, \int_{\la_{j-1}}^{\la_j} |f'|^2\,.
\]
Summing these estimates and recalling that $ \displaystyle \int_\R |f'|^2 =
4\pi^2\, \int_\R \xi^2 |\widehat{f}\,|^2 $, we get the claimed upper bound for the sum over $\la\in\Lambda$.
The proof of the upper bound for the sum over $\mu\in M$ is similar. \hfill $\Box$

\section{Uniqueness Pairs}
In this section we prove the following result, which somewhat extends the uniqueness part of
Theorem~\ref{thm-prelim}:

\medskip\noindent{\bf Theorem 1-UP.}
{\em Let $(\Lambda, M)$ be a supercritical pair with parameters $p$ and $q$.
Then the only function $f\in \Ho$ such that $f(\lambda) = \widehat{f}(\mu) = 0$, $\lambda\in\Lambda$, $\mu\in M$, is $f= 0$}.

\medskip
The proof splits into three  parts -- two parts that deal with smooth functions on $\R$,
and a complex-analytic part.
First, we prove that any function $f\in \Ho$ that vanishes on $\Lambda$ and whose Fourier transform vanishes on $M$ must  be a Schwartz function. Then we show that this function actually has an even faster decay together with its
Fourier transform and belongs to the Gelfand--Shilov space $\mathcal S(p, q)$, i.e.,
for some $c > 0$, $\displaystyle \int_\R |f(x)|^2e^{c|x|^p}\, {\rm d} x<\infty$ and $\displaystyle \int_\R |\widehat{f}(\xi)|^2e^{c|\xi|^q}\, {\rm d}\xi < \infty$.
In particular, $f$ and $\widehat{f}$ are entire functions of finite order.
The main tool in both steps is a modification of the Poincar\'e--Wirtinger inequality for fractional derivatives, which will be our starting point, and its slowly varying local version.

To finish the proof, we use some Phragm\'en--Lindel\"of type bounds from the entire function theory, working simultaneously on the space and the Fourier sides.

\subsection{From Sobolev to Schwartz spaces}\label{subsect:from Sobolev}
Our goal now is to prove that, under the assumptions of Theorem~1-UP,
for every $\theta > 1$, we have
$$\int_\R |f(x)|^2|x|^{2\theta}\, {\rm d}x < \infty \quad\text{and}\quad
\int_\R |\widehat{f}(\xi)|^2|\xi|^{2\theta}\, {\rm d}\xi < \infty$$
(recall that a function $f\in L^2(\R)$ belongs to Schwartz class $\mathcal S$
if and only if both of these integrals are finite for all $\theta>1$).

For convenience, we set $a = p - 1$, $b = q-1$,
where $p$ and $q$ are the parameters from the supercriticality assumption in
Theorem~1-UP. Then $ab=1$. Without loss of generality, we assume that $a \ge 1 \ge b$.

Next, observe that, by the supercriticality of the pair $(\Lambda, M)$,
there exist $\sigma > 1$ and $X_0, \Xi_0 > 0$ such that, for $x, \xi\in \R$,
the intervals of $\R \backslash \Lambda$ and $\R \backslash M$ containing $x$ and $\xi$
respectively have lengths at most $(2\sigma |x|^a)^{-1}$ and $(2\sigma |\xi|^b)^{-1}$, provided that
$|x| > X_0$, $|\xi| > \Xi_0$.
In order to apply Lemma~\ref{wirt2}, we take a non-negative
cutoff function $F\in C_0^\infty(-1, 1)$ with $\displaystyle \int_\R |F(x)|^2\, {\rm d}x = 1$,
and consider the function $g = fF_\upsilon$, where $F_\upsilon (x) = F(x-\upsilon)$, $\upsilon\in\R$.
If $|\upsilon| > X_0 + 1$, the function $g$ vanishes on a $(2\sigma)^{-1}(|\upsilon|-1)^{-a}$-dense set, so, applying Lemma~\ref{wirt2}(ii) to the function $g$, we get
\begin{equation}\label{temp2}
\Phi(\sigma^2(|\upsilon|-1)^{2a})\,
\int_\R |f(x)F_\upsilon(x)|^2\, {\rm d}x \le
\int_\R |(\widehat{f}*\widehat{F}_\upsilon)(\xi)|^2 \Phi(\xi^2)\, {\rm d}\xi
\end{equation}
with an arbitrary $C^1$-convex increasing function $\Phi$ vanishing at the origin.

Choosing $\Upsilon_0 > X_0 + 2$ so large that
$$ \left(\frac{\Upsilon_0+1}{\Upsilon_0 - 1}\right)^{2a} < \sigma, $$
and noting that $|x-\upsilon| \le 1$ for every $x$ in the support of $F_\upsilon$, we estimate the LHS of~\eqref{temp2} from below by
$$
\int_\R |f(x)|^2|F_\upsilon(x)|^2\Phi(\sigma|x|^{2a})\, {\rm d}x
\ge \sigma \int_\R |f(x)|^2 |F_\upsilon(x)|^2 \Phi( |x|^{2a})\, {\rm d}x,
$$
provided that $|\upsilon| > \Upsilon_0$. Plugging this into~\eqref{temp2} and
integrating both parts over such $\upsilon$, we obtain
\begin{equation}\label{temp3}
\sigma \int_{|x| > \Upsilon_0 + 1} |f(x)|^2 \Phi(|x|^{2a})\, {\rm d}x
\le \iint_{\R \times \R} |(\widehat{f}*\widehat{F_\upsilon})(\xi)|^2 \Phi(\xi^2)\,
{\rm d}\xi\, {\rm d}\upsilon.
\end{equation}

Now, we are going to bound the RHS from above. Recalling that $\widehat{F}_\upsilon(\eta) = \widehat{F}(\eta)e^{-2\pi {\rm i} \upsilon \eta}$ and using the Plancherel theorem, we have
\begin{equation*}
\int_\R \left|\int_\R \widehat{f}(\xi-\eta)\widehat{F}(\eta)e^{-2\pi {\rm i}\upsilon\eta}\,
{\rm d}\eta \right|^2\,{\rm d}\upsilon
= \int_\R |\widehat{f}(\xi-\eta)|^2|\widehat{F}(\eta)|^2\, {\rm d}\eta.
\end{equation*}
Thus, after a change of variables, the RHS of~\eqref{temp3} is equal to
\begin{equation}\label{temp12}
\iint_{\R\times\R}|\widehat{f}(\xi)|^2|\widehat{F}(\eta)|^2\Phi((\xi+\eta)^2)\,
{\rm d}\xi\, {\rm d}\eta.
\end{equation}

So far $\Phi$ was an arbitrary convex increasing $C^1$-function.
Now it is time to choose it.
We fix $\theta>1$ and a large number $T$, and put $\Phi_T(t) = t^\theta$, $0\le t \le T$, and
$\Phi_T(t) = T^\theta + \theta T^{\theta-1}(t-T)$, $t > T$.
Note that we have
\begin{equation}\label{temp4}
\Phi_T((\xi+\eta)^2) \le \Phi_T((|\xi|+|\eta|)^2) \le
\Phi_T(\xi^2) \left(1 + \frac{|\eta|}{|\xi|}\right)^{2\theta},
\end{equation}
which follows from the fact that $\Phi_T(t)t^{-\theta}$ is a non-increasing function.

Since $|\widehat{F}(\eta)|^2$ decays faster than any power of $\eta$,
for $|\xi|\to \infty$, we have
$$
\int_\R |\widehat{F}(\eta)|^2 \left( 1 + \frac{|\eta|}{|\xi|} \right)^{2\theta}\, {\rm d}\eta
\longrightarrow \int_\R |\widehat{F}(\eta)|^2\, {\rm d}\eta = \int_\R |F(x)|^2\, {\rm d}x = 1
$$
by the dominated convergence theorem, so we obtain
$$
\int_\R |\widehat{F}(\eta)|^2 \left(1 + \frac{|\eta|}{|\xi|}\right)^{2\theta}\, {\rm d}\eta < \sqrt{\sigma},
$$
as long as $|\xi| \ge \Xi_1 = \Xi_1(\sigma, p, \theta)$ (it is crucial that $\Xi_1$
does not depend on $T$). Thus,
\begin{equation}\label{temp5}
\iint_{\{|\xi|\ge \Xi_1\}\times \R}
|\widehat{f}(\xi)|^2|\widehat{F}(\eta)|^2\Phi_T((\xi+\eta)^2)\, {\rm d}\xi\, {\rm d}\eta
\le \sqrt{\sigma} \int_\R |\widehat{f}(\xi)|^2\Phi_T(\xi^2)\, {\rm d}\xi.
\end{equation}
Note that the integral on the RHS is finite since
$\displaystyle \int_\R |\widehat{f}(\xi)|^2(1+\xi^2)\, {\rm d}\xi~<~\infty$
(recall that, for every $T$, we have $\Phi_T(\xi^2) = O(\xi^2)$, as $|\xi|\to \infty$).

In~\eqref{temp12}, the integral over $\{|\xi|\le \Xi_1\}\times\R$ is dominated by
$$
\iint_{[-\Xi_1, \Xi_1]\times \R}
|\widehat{f}(\xi)|^2|\widehat{F}(\eta)|^2(\Xi_1 + |\eta|)^{2\theta}\, {\rm d}\xi\, {\rm d}\eta,
$$
which is some finite constant (independent of $T$) due to the fast decay of $\widehat{F}$.

At last, trivially, we have
$$
\int_{|x| \le \Upsilon_0 + 1}|f(x)|^2 \Phi_T(|x|^{2a})\,{\rm d}x \le
(\Upsilon_0 + 1)^{2a\theta}\, \int_\R |f(x)|^2\, {\rm d}x.
$$

Putting all the estimates together we get
\begin{equation}\label{eq-a}
\int_\R |f(x)|^2\Phi_T(|x|^{2a})\, {\rm d}x
\le \sigma^{-1/2}\int_\R |\widehat{f}(\xi)|^2\Phi_T(\xi^2)\, {\rm d}\xi + K,
\end{equation}
where $K = K_f > 0$ is independent of $T$.

\medskip
Now, let us switch the roles of $f$ and $\widehat{f}$ and consider the function
$\Psi_T(t) = \Phi_T(t^a)$ {(which is still convex, $C^1$-smooth, and increasing)}
instead of $\Phi_T$. Now, the function $\Psi_T(t)t^{-a\theta}$ is non-increasing,
the whole argument goes through, and we get
\begin{equation}\label{eq-b}
\int_\R |\widehat{f}(\xi)|^2\Psi_T(|\xi|^{2b})\, {\rm d}\xi
\le \sigma^{-1/2}\, \int_\R |f(x)|^2\Psi_T(x^2)\, {\rm d}x + K\,.
\end{equation}
The integral on the RHS of~\eqref{eq-b} coincides with the one on the LHS of~\eqref{eq-a},
and hence, is also finite.

\medskip
Now, we add inequalities~\eqref{eq-a} and~\eqref{eq-b}. Recalling that
$\Psi_T(|\xi|^{2b}) = \Phi_T(\xi^2)$ and $\Phi_T(|x|^{2a}) = \Psi_T(x^2)$ and
using the finiteness of the RHS, we obtain
$$
\int_\R |f(x)|^2\Phi_T(|x|^{2a})\, {\rm d}x +
\int_\R |\widehat{f}(\xi)|^2\Phi_T(\xi^2)\, {\rm d}\xi
\le \frac{2K}{1-\sigma^{-1/2}}\,,
$$
independently of the choice of $T>0$.
It remains to let $T\to \infty$ and use the monotone convergence theorem to get the bound
$$
\int_\R |f(x)|^2 |x|^{2a\theta}\, {\rm d}x +
\int_\R |\widehat{f}(\xi)|^2 |\xi|^{2\theta}\, {\rm d}\xi
\le \frac{2K}{1-\sigma^{-1/2}} < \infty.$$

\subsection{From Schwartz to Gelfand--Shilov spaces}\label{subsec:fromSchwartzToGSh}
Here, we show that if a Schwartz function $f$ vanishes on $\Lambda$ and its Fourier transform
$\widehat f$ vanishes on $M$, and the pair $(\Lambda, M)$ is supercritical with parameters $p$ and $q$, then, for every sufficiently big positive integer $\ell\ge \ell_0$, we have
\begin{equation}\label{eq-GS}
\int_\R |f(x)|^2 |x|^{p\ell}\, {\rm d}x, \
\int_\R |\widehat{f}(\xi)|^2 |\xi|^{q\ell}\, {\rm d}\xi
\le C^\ell \ell!\, \int_\R |f(x)|^2\, {\rm d}x\,.
\end{equation}
It is easy to see that these bounds yield convergence of the integrals
\footnote{Indeed, assume that $\| f \|_{L^2(\R)} = 1$.
Since the integrals $\displaystyle \int_\R |f(x)|^2 |x|^{p\ell}\, {\rm d}x$ are finite
for all $\ell\ge 0$, we can choose a constant $C_1$ so big that, for $0\le\ell <\ell_0$,
\[ \int_\R |f(x)|^2 |x|^{p\ell}\, {\rm d}x \le C_1^\ell \ell!\,.\]
Whence,
\[
\int_\R |f(x)|^2 e^{c|x|^p}\, {\rm d}x =
\sum_{\ell\ge 0} \frac{c^\ell}{\ell!}\, \int_\R |f(x)|^2 |x|^{p\ell}\, {\rm d}x
\le \sum_{\ell\ge 0} \bigl( c \cdot \max(C, C_1) \bigr)^\ell < \infty\,,
\]
provided that $c$ is chosen so small that $c \cdot \max(C, C_1) < 1$.
}
\[
\int_\R |f(x)|^2 e^{c|x|^p}\, {\rm d}x, \
\int_\R |\widehat{f}(\xi)|^2 e^{c|\xi|^q}\, {\rm d}\xi < \infty\,.
\]

In this part, we again rely on our version of the Poincar\'e--Wirtinger inequality
(Lemma~\ref{wirt2}). The main difference with the previous part is that now we need
to exercise some care in choosing the cutoff function $F$:
its choice has to be adjusted to the value of the parameter $\ell$ in~\eqref{eq-GS}.

\medskip
As in the first step, we set $a=p-1$, $b=q-1$, and
assume that $\Lambda$ is a discrete set such that for $|x| \ge X_0$ the constituting interval of $\R\backslash \Lambda$ containing $x$ has length at most
$(2\sigma |x|^a)^{-1}$ for some $\sigma > 1$.
Let $f\in \So$ vanish on $\Lambda$. Fix $\theta\ge 1$ and consider the function
$$
F = \done_{[-u, u]}*\underbrace{\tfrac{k}{u}\done_{[-\frac{u}{2k}, \frac{u}{2k}]}* \ldots  * \tfrac{k}{u}\done_{[-\frac{u}{2k}, \frac{u}{2k}]}}_{k\text{--fold\ convolution}}\,,
$$
where $\done_I$ is the indicator function of the interval $I$, and large
parameters $k$ and $u$ are to be chosen later.
We are going to apply Lemma~\ref{wirt2}(ii) with $\Phi(x) = x^\theta$
to the function $fF_\upsilon$, where $F_\upsilon(x) = F(x-\upsilon)$.
Note that $\operatorname{supp}(F_\upsilon) \subset [\upsilon - \frac32 u, \upsilon + \frac32 u]$.
If $|\upsilon|-\frac{3}{2}u \ge X_0$, the function $fF_\upsilon$ vanishes on a
$(2\sigma)^{-1}\, (|\upsilon|-\frac{3}{2}u)^{-a}$-dense set, so Lemma~\ref{wirt2}(ii) applied with
$t = \sigma (|\upsilon|- \frac32 u)^{2a}$ gives us
\begin{equation}\label{temp13}
\Bigl( |\upsilon|-\frac{3}{2}u \Bigr)^{2a\theta}\, \int_\R |f(x)|^2|F_\upsilon(x)|^2\, {\rm d}x
\le \sigma^{-2\theta}\, \int_\R |\xi|^{2\theta}|(\widehat{f}*\widehat{F}_\upsilon)(\xi)|^2\, {\rm d}\xi.
\end{equation}

Now, we integrate both sides over $|\upsilon| \ge \tfrac{3}{2}u + X_0$. The integral of the LHS
is at least
$$
\int_{|x|\ge 3u + X_0}(|x|-3u)^{2a\theta}|f(x)|^2\, {\rm d}x
\cdot \int_\R |F(x)|^2\, {\rm d}x,
$$
because, as long as $x$ with $|x|\ge 3u + X_0$ is contained in $\supp F_\upsilon$, we have
$|\upsilon|\ge \tfrac{3}{2}\,u+X_0$ and $|\upsilon|-\tfrac{3}{2}\,u \ge |x|-3u$.
Assuming $u \ge X_0$ and $K > 4$, we can bound this quantity from below by
$$
\Bigl( \frac{K-3}{K} \Bigr)^{2a\theta}\, \int_{|x|\ge Ku} |x|^{2a\theta}|f(x)|^2\, {\rm d}x
\cdot \int_\R |F(x)|^2\, {\rm d}x.
$$

Let us now estimate the integral
$\displaystyle
\int_\R{\rm d}\upsilon\, \int_\R |\xi|^{2\theta}|(\widehat{f}*\widehat{F_\upsilon})(\xi)|^2\,
{\rm d}\xi $,
which appears when we integrate the RHS of~\eqref{temp13} over $\upsilon$.
As in the previous step, noting that
$\widehat{F_\upsilon}(\eta) = e^{-2\pi {\rm i} \upsilon \eta}\widehat{F}(\eta)$
and using the Plancherel theorem, we find that it is equal to
$$
\iint_{\R\times\R}
|\widehat{f}(\xi)|^2 |\widehat{F}(\eta)|^2 |\xi+\eta|^{2\theta}\, {\rm d}\xi\, {\rm d}\eta\,.
$$
Note that
$$
\widehat{F}(\eta) = \widehat{\done}_{[-u, u]}(\eta) \cdot
\Bigl(\, \frac{\sin \pi \frac{u}{k}\eta}{\pi \frac{u}{k}\eta} \,\Bigr)^k.
$$
We claim that
$$
(|\xi|+|\eta|) \cdot \Bigl|
\frac{\sin \pi \frac{u}{k}\eta}{\pi \frac{u}{k}\eta} \,\Bigr| \le |\xi| + \frac{k}{\pi u}.$$
If $|\eta| \le k/(\pi u)$ the claim becomes obvious after we estimate the second factor on the LHS by $1$. Otherwise, we estimate the absolute value of the sine on the LHS by $1$, and note that
$$
\pi \frac{u|\eta|}{k} \cdot \Bigl( |\xi| + \frac{k}{\pi u}\, \Bigr) \ge |\xi| + |\eta|
$$
because in this case $\pi u |\eta|/k \ge 1$. Thus, if $k\ge \theta$, we have
\begin{align*}
\int_\R \int_\R |\widehat{f}(\xi)|^2 |\widehat{F}(\eta)|^2|\xi &+\eta|^{2\theta}\,
{\rm d}\xi\, {\rm d}\eta \\
& \le \int_\R \int_\R |\widehat{f}(\xi)|^2 \left(|\xi|+\frac{k}{\pi u}\right)^{2\theta}
|\widehat{\done}_{[-u, u]}(\eta)|^2\, {\rm d}\xi \, {\rm d}\eta \\
&=2u\int_\R |\widehat{f}(\xi)|^2\left(|\xi|+\frac{k}{\pi u}\right)^{2\theta}\, {\rm d}\xi,
\end{align*}
where in the first step we again used that $|\sin t|/|t| \le 1$.

Taking into account that $\displaystyle \int_\R |F(x)|^2\, {\rm d}x \ge u$ (just because $F(x) = 1$ if  $|x| \le u/2$) and juxtaposing the estimates for the LHS and RHS of~\eqref{temp13}, we arrive at the inequality
\begin{equation}\label{eq:temp50}
\int_{|x| \ge Ku} |x|^{2a\theta}|f(x)|^2\, {\rm d}x
\le 2\Bigl(\, \frac1{\sigma}\, \Bigl(\, \frac{K}{K-3}\, \Bigr)^a \Bigr)^{2\theta}\,
\int_\R |\widehat{f}(\xi)|^2 \Bigl(|\xi|+\frac{k}{\pi u}\, \Bigr)^{2\theta}\, {\rm d}\xi,
\end{equation}
provided that $u\ge X_0$, $K > 4$, and $k\ge \theta$.
We will be interested in large $\theta$. So, we choose $K$ so large that
$\frac1\sigma\, \bigl(\, \frac{K}{K-3}\,\bigr)^a  < 1$.
Then we shall take $k\in [\theta, \pi \theta]$ and split the integral on the RHS into the sum
of the integrals
\begin{align*}
\Bigl( \int_{|\xi| > K \frac{k}{\pi u}} &+ \int_{|\xi| \le K \frac{k}{\pi u}}  \Bigr)
|\widehat{f}(\xi)|^2 \Bigl(|\xi|+\frac{k}{\pi u}\, \Bigr)^{2\theta}\, {\rm d}\xi \\
&\le
\Bigl( \frac{K+1}K \Bigr)^{2\theta}\,
\int_{|\xi| \ge K \frac{k}{\pi u}} |\widehat{f}(\xi)|^2 |\xi|^{2\theta}\, {\rm d}\xi
+ \int_{|\xi| \le K \frac{k}{\pi u}}
|\widehat{f}(\xi)|^2 \Bigl(|\xi|+\frac{k}{\pi u}\, \Bigr)^{2\theta}\, {\rm d}\xi
\\
&\le \Bigl( \frac{K+1}K \Bigr)^{2\theta}\,
\int_\R |\widehat{f}(\xi)|^2 |\xi|^{2\theta}\, {\rm d}\xi
+ \int_{|\xi| \le K \frac{k}{\pi u}}
|\widehat{f}(\xi)|^2 \Bigl(|\xi|+\frac{k}{\pi u}\, \Bigr)^{2\theta}\, {\rm d}\xi\,.
\end{align*}
Adding the integral $\displaystyle \int_{|x| \le Ku} |x|^{2a\theta} |f(x)|^2\, {\rm d}x$
to both sides of~\eqref{eq:temp50}, we obtain
\begin{multline*}
\int_\R |x|^{2a\theta}|f(x)|^2\, {\rm d}x \le
2\left( \tfrac1\sigma\, \tfrac{K+1}{K}\left(\tfrac{K}{K-3}\right)^a\right)^{2\theta}\,
\int_\R |\widehat{f}(\xi)|^2|\xi|^{2\theta}\, {\rm d}\xi
\\
+ 2\int_{|\xi|\le K\frac{k}{\pi u}}
|\widehat{f}(\xi)|^2 \Bigl(|\xi|+\frac{k}{\pi u}\Bigr)^{2\theta}\, {\rm d}\xi
+ \int_{|x| \le Ku} |x|^{2a\theta}|f(x)|^2\, {\rm d}x\,.
\end{multline*}
Bounding the second and the third terms on the RHS as
\[
\int_{|\xi|\le K\frac{k}{\pi u}}
|\widehat{f}(\xi)|^2 \Bigl(|\xi|+\frac{k}{\pi u}\Bigr)^{2\theta}\, {\rm d}\xi
\stackrel{k\le \pi\theta}\le (K+1)^{2\theta} \Bigl(\frac{\theta}{u}\Bigr)^{2\theta}\,
\int_\R |f|^2
\]
and
\[
\int_{|x| \le Ku} |x|^{2a\theta}|f(x)|^2\, {\rm d}x
\le K^{2a\theta}u^{2a\theta}\, \int_\R |f|^2\,,
\]
and assuming further that $\frac{K+1}{K}\left(\frac{K}{K-3}\right)^a\sigma^{-1} < 1$, we conclude that,
for large enough $\theta$, we have
\begin{multline*}
\int_\R |x|^{2a\theta}|f(x)|^2\, {\rm d}x \le
\frac{1}{2}\int_\R |\widehat{f}(\xi)|^2|\xi|^{2\theta}\, {\rm d}\xi
\\
+\Bigl(K^{2a\theta}u^{2a\theta} + 2(K+1)^{2\theta}\Bigl(\frac{\theta}{u}\Bigr)^{2\theta}\Bigr)\,
\int_\R |f(x)|^2\, {\rm d}x.
\end{multline*}
At last, choosing $u = \theta^{\tfrac{1}{a+1}}$, which is consistent with $u \ge X_0$ as
$\theta\to \infty$, we finally get
\begin{equation}\label{temp14}
\int_\R |x|^{2a\theta}|f(x)|^2\, {\rm d}x
\le \frac{1}{2}\, \int_\R |\xi|^{2\theta} |\widehat{f}(\xi)|^2\, {\rm d}\xi
+ C^{2\theta}\theta^{\tfrac{2a}{a+1}\theta}\, \int_\R |f(x)|^2\, {\rm d}x.
\end{equation}

By the very same argument if $\widehat{f}$ vanishes on the set $M$ such that, for $|\xi|\ge \Xi_0$, the constituting interval of $\R\backslash M$ containing $\xi$ has length at most $(2\sigma |\xi|^b)^{-1}$ with $\sigma > 1$ and $b = a^{-1}$, we obtain
\begin{equation}\label{temp15}
\int_\R |\xi|^{2b\kappa}|\widehat{f}(\xi)|^2\, {\rm d}\xi
\le \frac{1}{2}\,\int_\R |x|^{2\kappa}|f(x)|^2\, {\rm d}x
+ C^{2\kappa}\kappa^{\tfrac{2b}{b+1}\kappa}\, \int_\R |\widehat{f}(\xi)|^2\, {\rm d}\xi,
\end{equation}
provided that $\kappa$ is large enough.

We will choose $\theta$ and $\kappa$ so that $\kappa = a \theta$, and, since $b=a^{-1}$,
$\theta = b\kappa$. Note that in this case we have
\[
\frac{a}{a+1}\, \theta = \frac{\kappa}{a+1} = \frac{b}{b+1}\, \kappa,
\]
so the second terms in the RHS of~\eqref{temp14} and~\eqref{temp15} are the same up to the value of
the constant $C$. Adding~\eqref{temp14} and~\eqref{temp15},
we thus obtain, for big enough $\theta$,
$$
\int_\R |x|^{2a\theta}|f(x)|^2\, {\rm d}x
+ \int_\R |\xi|^{2b\kappa}|\widehat{f}(\xi)|^2\, {\rm d}\xi
\le C_1^{2\theta}\theta^{\tfrac{2a}{a+1}\theta}\, \int_\R |f(x)|^2\, {\rm d}x.
$$
At last, we choose
\[
\theta = \frac{a+1}{2a}\, \ell, \quad \ell = \ell_0, \ell_0+1, \ell_0+2, \ldots \,,
\]
where $\ell_0$ is a sufficiently large positive integer.
Since $a=p-1$, with this choice of the parameters, we have
$2a\theta = p\ell$, $2\theta = q\ell$, and therefore, for $\ell \ge \ell_0$,
\[
\max\Bigl[ \int_\R |x|^{p\ell}|f(x)|^2\, {\rm d}x,
\int_\R |\xi|^{q\ell}|\widehat{f}(\xi)|^2\, {\rm d}\xi \Bigr]
\le C_2^\ell \ell^\ell\,
\int_\R |f(x)|^2\, {\rm d}x \le C_3^\ell \ell!\, \int_\R |f(x)|^2\, {\rm d}x\,,
\]
proving~\eqref{eq-GS}.

\subsection{Phragm\'en--Lindel\"of tricks} We know that, for some $c > 0$, the integrals
$\displaystyle \int_\R |f(x)|^2 e^{c|x|^p}\, {\rm d}x $ and
$ \displaystyle \int_\R |\widehat{f}(\xi)|^2 e^{c|\xi|^q}\, {\rm d}\xi $
are finite.
These estimates imply, in particular, that $f'$ and $\widehat{f}'$ are bounded in absolute value by some
constant $A$. Therefore, for every $x_0\in\R$,
\begin{align*}
\int_\R |f(x)|^2 e^{c|x|^p}\, {\rm d}x &\ge
\int_\R \max(|f(x_0)| - A|x-x_0|, 0)^2 e^{c|x|^p}\, {\rm d}x
\\
&\ge
\int_{x_0 - (2A)^{-1}|f(x_0)|}^{x_0 + (2A)^{-1}|f(x_0)|} (|f(x_0)| - A(x-x_0))^2 e^{c|x|^p}\, {\rm d}x
\ge \frac1{8A}\, |f(x_0)|^3 e^{c|x_0|^p}\,,
\end{align*}
whence,
$|f(x)| = O(e^{-\frac{c}3 |x|^p})$ as $x\to\infty$ and, similarly,
$|\widehat{f}(y)| = O(e^{-\frac{c}3 |y|^q})$ as $y\to\infty$.
We also note that, by standard estimates of the Fourier transform,
$f$ is an entire function of order $p$, and $\widehat{f}$ is an entire function of order $q$, and naturally, in this part we will use some complex analysis techniques.

\smallskip
We begin by defining the pair of the Phragm{\'e}n--Lindel{\"o}f indicators
$$k_1(\theta) = \frac{1}{2\pi}\, \limsup_{r\to \infty} \frac{\log |f(re^{i\theta})|}{r^p}\,,$$
and
$$k_2(\phi) = \frac{1}{2\pi}\, \limsup_{\pho\to\infty} \frac{\log |\widehat{f}(\pho e^{i\phi})|}{\pho^q}$$
(the factors $(2\pi)^{-1}$ are introduced to simplify some expressions).
We note that both indicators are negative at $0$ and $\pm\pi$, and set
$\kappa_i = \min (|k_i(0)|, |k_i(\pm \pi)|)$, $i = 1, 2$.

\begin{claim}\label{claim1}
We have
$$k_1(\theta) \le \frac{1}{p(q\kappa_2)^{p/q}}|\sin \theta|^p,$$
and
$$k_2(\phi) \le \frac{1}{q(p\kappa_1)^{q/p}}|\sin \phi|^q.$$
\end{claim}

\noindent{\em Proof:}
For every $\eps > 0$, we have
\begin{align*}
|f(z)| &\le \int_\R |\widehat{f}(\xi)|e^{2\pi |\xi|\, |{\rm Im}z|}\, {\rm d}\xi \\
&\le C_\eps\, \int_\R e^{(-2\pi \kappa_2+2\eps)|\xi|^q + 2\pi |\xi|\, |{\rm Im}z|}
e^{-\eps |\xi|^q}\, {\rm d} \xi \\
&\le C_{\eps}'e^{2\pi \max\limits_{\xi\in \R} (|\xi|\, |{\rm Im}z| - (\kappa_2-\frac{\eps}{\pi})|\xi|^q)}\, .
\end{align*}
By the Young inequality,
\[
|\xi|\, |{\rm Im}z| \le \bigl( \kappa_2-\frac{\eps}{\pi} \bigr)|\xi|^q +
\frac{1}{p}\, \frac{|{\rm Im}z|^p}{( q(\kappa_2-\frac{\eps}{\pi}))^{p/q}}\,.
\]
Hence, setting $z=re^{{\rm i}\theta}$ and letting $r\to\infty$, we get
$$k_1(\theta) \le \frac{1}{p}\, \frac{|\sin \theta|^p}{(q(\kappa_2-\frac{\eps}{\pi}))^{p/q}}\,.$$
Since $\eps >0$ was arbitrary, this gives the first bound. The proof of the second estimate is the same.
\hfill $\Box$

\medskip
Replacing $f(x)$ by $f(-x)$ if necessary, we may assume without loss of generality that $\kappa_1 = |k_1(0)|$.
Replacing $f(x)$ by $\overline{f(x)}$ after that, if needed, we can also assume that $\kappa_2 = |k_2(0)|$
as well. We also note that both quantities $\kappa_1$, $\kappa_2$
are finite (the Phragm\'en--Lindel\"of principle
yields that if the indicator of an entire function of finite order and bounded type
equals to $-\infty$ at one point, then the function vanishes identically).

To use the density of zeroes of $f$ and $\widehat f$, we need an observation about
analytic functions of exponential type in an angle $\Gamma(\alpha) = \{z\colon |{\rm arg}(z)|\le \alpha\}$ with $\alpha< \pi/2$.

\begin{claim}\label{claim-observation}
Suppose that $F$ is analytic in $\Gamma (\alpha)$, continuous up to the boundary, and satisfies
$|F(z)|\le Me^{a\pi |{\rm Im}z|}$ with some $a<1$. Let $\mathfrak Z = \{\mathfrak z_1, \mathfrak z_2, \dots \} \subset \R_+$ be a discrete set with $\mathfrak z_1 < \mathfrak z_2 < \ldots $ and
$\mathfrak z_{j+1}-\mathfrak z_j \le 1$. If $F|_\mathfrak Z = 0$, then
\[
\limsup_{x\to +\infty} \frac{\log |F(x)|}x \le Q(a, \alpha)<0.
\]
\end{claim}
\noindent Note that we do not require $\mathfrak z_1$ to be close to the origin and that
it does not appear as an argument of $Q$.

\medskip\noindent{\em Proof:}
Take any point $x\in \R_+$ and consider the disk $D$ of radius $r=x\sin\alpha$ centered at $x$. Note that $D\subset \Gamma(\alpha)$. We shall take $x$ so large that
$(1-\sin\alpha)x>\mathfrak z_1$. Then the points of the set $\mathfrak Z \cap D$ form a $1$-net
on the diameter of $D$. By the classical Jensen formula, we have
\[
\log |F(x)| \le \frac1{2\pi}\, \int_{-\pi}^\pi \log|F(x+re^{{\rm i}\eta})|\, {\rm d}\eta
+ \sum_{\mathfrak z \in \mathfrak Z \cap D} \log\frac{|x-\mathfrak z|}r\,.
\]
The first term on the RHS is bounded from above by
\[
\frac1{2\pi}\, \int_{-\pi}^\pi a\pi r|\sin\eta|\, {\rm d}\eta + \log M = 2ar + O(1),
\]
while the sum can be estimated from above by
\[
2 \sum_{1\le k \le r} \log\frac{k}r
\le 2 \sum_{1\le k \le r} \int_k^{k+1} \log\frac{t}r\, {\rm d}t
\le 2 \int_1^{r+1} \log\frac{t}r\, {\rm d}t
= (-2 + o(1)) r
\]
as $r\to\infty$. Thus,
\[
\frac1x\, \log |F(x)| \le (2(a-1)+o(1))\sin\alpha\,, \quad x\to\infty\,,
\]
which implies the claim with $Q(a, \alpha) = -2(1-a)\sin\alpha$.
\hfill $\Box$

\medskip
Note that if elements of $\mathfrak Z$ are $\ell$-dense instead of being $1$-dense, we can consider
$F(\ell z)$ instead of $F$, assume that $a\,\ell <1$, and arrive at essentially the same conclusion.

\medskip
Now we claim that, for all $\theta\in \displaystyle \bigl( 0, \frac{\pi}{2p} \bigr)$, we have
\begin{equation}\label{eq:ineq}
\kappa_1 \ge \frac1p\, \tan p\theta - \frac1{p(q\kappa_2)^{p/q}}\,
\frac{(\sin\theta)^p}{\cos p\theta}\,.
\end{equation}
Indeed, consider the function $F(z) = e^{(2\pi \kappa_1 - \eps)z} f(z^{1/p})$ with
small $\eps>0$ in the right half-plane ${\rm Re}(z) >0$. By Claim~\ref{claim1},
for $|\theta| < \displaystyle \frac{\pi}{2p}$, we have
\begin{multline*}
\log |F(re^{{\rm i}p\theta})| \le
2\pi \bigl( \kappa_1 \cos p\theta + k_1(\theta) \bigr)r + O(1) \\
\le 2\pi \bigl( \kappa_1 \cos p\theta + \frac1{p(q\kappa_2)^{p/q}}\, |\sin\theta|^p
\bigr)\, r + O(1)\,.
\end{multline*}

Suppose that estimate~\eqref{eq:ineq} does not hold for some $\theta_0\in
\displaystyle \bigl( 0, \frac{\pi}{2p} \bigr)$, that is,
\[
\kappa_1 \cos p\theta_0 + \frac1{p(q\kappa_2)^{p/q}}\, (\sin\theta_0)^p
< \frac{c}p\, \sin p\theta_0
\]
with some $c<1$.
Consider the function $F$ in the angle $\Gamma (\alpha)$ with $\alpha = p\theta_0$.
On the boundary rays $|F(z)| = O\bigl( e^{\frac{2\pi c}p\,|z|\sin p\theta_0}\bigr)
= O(e^{a\pi |{\rm Im}z|})$ with $a=\frac{2c}p$ satisfying $a\, \frac{p}2<1$.
On the positive ray $|F|$ remains bounded (recall that $\kappa_1 = - k_1(0)$). At last,
$|F(z)| = O(e^{C|z|})$ everywhere in $\Gamma (\alpha)$. Therefore,
by the Phragm\'en--Lindel\"of principle,
the function $F(z)e^{a\pi {\rm i}z}$ is bounded in the angle $\{0 \le {\rm arg}(z) \le \alpha\}$
(the upper half of $\Gamma(\alpha)$). Thus, we have $\log |F(z)| \le a\pi {\rm Im} z + O(1)$
for ${\rm Im}z \ge 0$, $z\in\Gamma (\alpha)$. Arguing similarly in the lower half of $\Gamma (\alpha)$,
we see that $|F(z)| \le O(e^{a\pi |{\rm Im} z|})$ everywhere in $\Gamma (\alpha)$.

On the other hand, $F$ vanishes at every point $\la^p$, $\la\in\Lambda\cap (0, \infty)$. Since
$\la_{j+1}^p - \la_j^p \le p\la_{j+1}^{p-1}(\la_{j+1}-\la_j) < \tfrac12 p$ for large enough $j$ by the supercriticality of $\Lambda$, applying Claim~\ref{claim-observation}, we conclude that
\[
\limsup_{x\to +\infty} \frac{\log |F(x)|}x \le Q < 0
\]
{\em independently of}  $\eps$. However, the definition of $\kappa_1 = - k_1(0)$ implies that it is impossible if $\eps < |Q|$, proving estimate~\eqref{eq:ineq}.

\medskip
Thus, we have
$$
\kappa_1 \ge \frac{1}{p}\tan p\theta - \frac{1}{(q \kappa_2 )^{p/q}}\,
\frac{|\sin \theta|^p}{p \cos p\theta}\,, \quad 0 < \theta < \frac{\pi}{2p}\,.
$$
Next, we will need the trigonometric inequality
$$\frac{1}{p}\tan p\theta \ge \frac{\sin \theta}{(\cos p\theta)^{1/p}}, \quad
0 \le \theta < \frac{\pi}{2p}\,,$$
with strict inequality for $\theta > 0$. Indeed, both sides vanish at $\theta = 0$ and the derivative of the LHS is $\displaystyle \frac{1}{\cos^2 p\theta}$, while the derivative of the RHS is
$$\frac{\cos \theta \cos p\theta + \sin \theta \sin p\theta}{(\cos p\theta)^{1+1/p}} = \frac{\cos (p-1)\theta}{(\cos p\theta)^{1+1/p}}\,,$$
which is clearly smaller for $\theta\in (0, \tfrac{\pi}{2p})$.

Letting $\displaystyle \eta = \frac{\sin \theta}{(\cos p\theta)^{1/p}}$, we get
$$\kappa_1 > \eta - \frac{1}{(q \kappa_2)^{p/q}}\, \frac{\eta^p}{p}\,.$$
Maximizing the RHS with respect to $\eta$ (note that $\eta$ runs through the whole
interval $(0, +\infty)$ when $\theta$ increases from $0$ to $\pi/(2p)$),
we let $\eta = q \kappa_2$, which ultimately gives us $\kappa_1 > \kappa_2$.

The same argument with the roles of $p$ and $q$ exchanged gives us $\kappa_2 > \kappa_1$.
The contradiction we arrive at shows that $f$ is a zero function, completing the proof of Theorem~1-UP. \hfill $\Box$

\section{Frame bounds}
In this Section we turn to Theorem~2, using ideas
similar to the ones we used in the proofs of Theorem~2A
dealing with uniformly supercritical sets and in the first step of the proof of Theorem~1-UP.
{\em Throughout this section we always assume that $(\Lambda, M)$ is a separated supercritical
pair with parameters $(p, q)$, $\tfrac1p + \tfrac1q = 1$}.

First, we show that the lower frame-type bound in the conclusion of Theorem~\ref{thm-stable-uniq} holds
if we add to each of the sets $\Lambda$, $M$ finitely many points. Then, using a classical Duffin--Schaeffer argument, we remove this extra assumption and prove this lower bound without it. Finally, we prove the upper frame-type bound and thus deduce Theorem~\ref{thm-stable-uniq} in full generality.

\subsection{Adding finitely many points to $\Lambda$ and $M$}\label{subsection:Adding}
Here, we show that for all $s > 0$ such that $s\min (p, q) \ge 1$ there exist sets $\Lambda'\supset \Lambda$,
$M'\supset M$ and a constant $C > 0$, such that $\Lambda'\backslash \Lambda$ and $M'\backslash M$ are finite, and, for  all functions $f\in \Ho$,
we have
\begin{equation}\label{eq:*}
||f||_{\Ho_{s, p, q}}^2 \le C \Bigl(
\sum_{\lambda\in \Lambda'} (1+|\lambda|)^{(2s-1)p+1} |f(\lambda)|^2
+ \sum_{\mu \in M'} (1+|\mu|)^{(2s-1)q+1} |\widehat{f}(\mu)|^2 \Bigr).
\end{equation}

\medskip\noindent{\em Proof of~\eqref{eq:*}}:
We assume that both sums on the RHS are finite (otherwise, there is nothing
to prove).
We will use the same notation and conventions as in Section~\ref{subsect:from Sobolev}. Specifically, we assume without loss of generality that $p \ge 2 \ge q$ and let $a = p-1$, $b = q-1$.
We consider the same convex increasing functions $\Phi_T$ and $\Psi_T$ as in that section:
for big number $T$, we put $\Phi_T(t) = t^\theta, 0 \le t \le T$ and $\Phi_T(t) = T^\theta + \theta T^{\theta-1}(t-T)$, $t\ge T$, where $\theta = qs$, and $\Psi_T(t) = \Phi_T(t^a)$. Note that for  $t \ge 0$ we have $\Phi_T(t)\le t^\theta$, $\Phi_T'(t)\le \theta t^{\theta-1}$, and $\Psi_T(t) \le t^{a\theta}$, $\Psi_T'(t)\le a\theta t^{a\theta - 1}$.

Let $F\in C_0^\infty(-1, 1)$ with $\displaystyle \int_\R |F(x)|^2\, {\rm d}x = 1$. Put
$F_\upsilon(x)= F(x-\upsilon)$.
For big $x$, the constituting interval of $\R\backslash \Lambda$ containing $x$ has length at most $(1-\eps)(2\sigma (1+|x|)^{p-1})^{-1}$ for some $\eps > 0$, $\sigma > 1$. The set $\Lambda'$ will be $\Lambda$ with a lot of points added  near the origin so that all the constituting intervals of $\R\backslash\Lambda'$ have length bounded by a small numerical constant. Applying Lemma~\ref{wirt2}(i)
to the function $fF_\upsilon$ and the function $\Phi_T$, we get
\begin{multline}\label{eq:new}
\Phi_T(\sigma^2U(\upsilon)^2)\int_\R |f(x)|^2|F_\upsilon(x)|^2\, {\rm d}x \le \,
\int_\R |(\widehat{f}*\widehat{F}_\upsilon)(\xi)|^2\Phi_T(\xi^2)\, {\rm d}\xi \\
+ C_\eps \sigma U(\upsilon)\Phi_T'(\sigma^2U(\upsilon)^2)\sum_{\lambda\in \Lambda'} |f(\lambda)|^2|F_\upsilon(\lambda)|^2,
\end{multline}
where $U(\upsilon) =\max(K, |\upsilon|)^{p-1}$ and positive $K$ can be chosen arbitrary large,
provided that we add sufficiently many points to $\Lambda$.

To proceed, we need an auxiliary estimate.
\begin{claim}\label{claim:aux}
For $|x-\upsilon|\le 1$, we have
\[
\Phi_T(\sigma^2U(\upsilon)^{2})\ge (1-\eps_1)\Phi_T(\sigma^2U(x)^{2}),
\]
where
$\eps_1$ does not depend on $T$ and can be chosen arbitrarily small when $K$ is big enough,
while the parameters $p$, $q$ and $s$ remain fixed.
\end{claim}

\noindent{\em Proof of Claim~\ref{claim:aux}}:
Since the function $U$ is even, we may assume that $x$ and $\upsilon$ are non-negative.
If $\upsilon\ge x$ or $\upsilon \le x \le K$ the estimate is obvious with $\eps_1=0$.
In the remaining case, we have $K-1<x-1\le \upsilon<x$, whence
\[
U(\upsilon) \ge U(x-1) \ge \bigl( x-1)^{p-1} = (1-x^{-1})^{p-1} U(x) \ge (1-K^{-1})^{p-1} U(x).
\]
Next, we note that $\Phi_T'(t) \le \theta t^{-1} \Phi_T(t)$, whence (taking into account that
$\Phi_T'$ does not decrease), for $t'<t$,
\[
\Phi_T(t)\le \Phi_T(t') + (t-t') \Phi_T(t)
\le \Phi_T'(t') + \theta \bigl( 1-\frac{t'}{t} \bigr) \Phi_T(t).
\]
Using this with $t=\sigma^2  U(x)^2$ and $t'=\sigma^2 U(\upsilon)^2$ (and recalling that
$\theta=qs$ is fixed), we get what we wanted. \hfill $\Box$

\medskip
Now, integrating inequality~\eqref{eq:new} over $\upsilon\in \R$, we obtain
\begin{multline}\label{eq:new-new}
\iint_{\R\times\R}
|f(x)|^2|F_\upsilon(x)|^2\Phi_T(\sigma^2U(\upsilon)^{2})\, {\rm d}\upsilon\, {\rm d}x
\le \iint_{\R\times\R}
|(\widehat{f}*\widehat{F}_\upsilon)(\xi)|^2\Phi_T(\xi^2)\, {\rm d}\xi\, {\rm d}\upsilon \\
+ C_\eps \sigma \sum_{\lambda\in \Lambda'} |f(\lambda)|^2
\int_{|\upsilon-\lambda|\le 1}
U(\upsilon)\Phi_T'(\sigma^2U(\upsilon)^2) |F_\upsilon(\lambda)|^2\, {\rm d}\upsilon.
\end{multline}
For the LHS of~\eqref{eq:new-new}, using Claim~\ref{claim:aux}, we have
$$
\iint_{\R\times\R}
|f(x)|^2|F_\upsilon(x)|^2\Phi_T(\sigma^2U(\upsilon)^{2})\, {\rm d}\upsilon\, {\rm d}x
\ge (1-\eps_1)\, \int_\R |f(x)|^2\Phi_T(\sigma^2U(x)^{2})\, {\rm d}x.
$$
Since $\Phi_T$ is convex and $\Phi_T(0) = 0$, we have
$$(1-\eps_1)\, \int_\R |f(x)|^2\Phi_T(\sigma^2U(x)^{2})\, {\rm d}x \ge (1-\eps_1)\sigma^2\,
\int_\R |f(x)|^2\Phi_T(U(x)^{2})\, {\rm d}x\,.$$

Using that $\Phi_T'(\sigma^2U(\upsilon)^2) \le \theta \sigma^{2\theta-2}U(\upsilon)^{2\theta - 2}$ and that, for $|\lambda-\upsilon|\le 1$, $U(\upsilon)^{2\theta-1}\le C_1 (K) U(\lambda)^{2\theta-1}$,
we bound from above the sum term which appears on the RHS of~\eqref{eq:new-new} by
$$
C_1(K) \cdot C_\eps \theta \sigma^{2\theta-1}\, \sum_{\lambda\in \Lambda'}|f(\lambda)|^2U(\lambda)^{2\theta-1}.
$$

As in Sections~\ref{subsect:from Sobolev} and~\ref{subsec:fromSchwartzToGSh}, by the Plancherel theorem, after integration over $\upsilon$
the integral term on the RHS of ~\eqref{eq:new-new} equals
\begin{equation}\label{eq:integral}
\iint_{\R\times \R} |\widehat{f}(\xi)|^2 |\widehat{F}(\eta)|^2\Phi_T(|\xi+\eta|^2)\,
{\rm d}\xi\,{\rm d}\eta.
\end{equation}
To bound this integral, we need another auxiliary estimate for the function $\Phi_T$.
\begin{claim}\label{claim:auxPhi}
For all $\delta > 0$, there exists $C_{\delta}$ (possibly also depending on $s, p, q$ but not on $T$) such that for all $\xi, \eta\in \R$,
we have
\begin{equation}\label{Phi bound}
\Phi_T(|\xi+\eta|^2)\le
(1+\delta)\Phi(\xi^2) + C_\delta |\eta|^{2\theta}.
\end{equation}
\end{claim}

\noindent{\em Proof of Claim~\ref{claim:auxPhi}}:
We use the bound~\eqref{temp4}
\[
\Phi_T((\xi+\eta)^2)
\le \Phi_T(\xi^2) \left(1 + \frac{|\eta|}{|\xi|}\right)^{2\theta}\,.
\]
If $|\eta|\le c |\xi|$, this estimate gives us
$\Phi_T((\xi+\eta)^2) \le \Phi_T(\xi^2) (1+c)^{2\theta} < (1+\delta)\Phi_T(\xi^2)$,
provided that $(1+c)^{2\theta} < 1+\delta$.
If $|\eta|\ge c|\xi|$, then $\Phi_T ((\xi+\eta)^2) \le \Phi_T ((1+c^{-1})^2 \eta^2)
\le (1+c^{-1})^{2\theta} |\eta|^{2\theta}$. \hfill $\Box$

\medskip
Now, we return to the integral~\eqref{eq:integral}.
Since $|\widehat{F}(\eta)|^2$ decays faster than any negative power and
$\displaystyle \int_\R |\widehat{F}(\eta)|^2\, {\rm d}\eta = \int_\R |F(x)|^2\, {\rm d}x = 1$,
using~\eqref{Phi bound} and letting $V(\xi) = \max(K, |\xi|)^{q-1}$,  we get
\begin{multline*}
\int_\R |\widehat{F}(\eta)|^2\Phi_T(|\xi+\eta|^2) {\rm d}\eta \le (1+\delta)\Phi_T(\xi^2) + C_{\delta}' \\
= (1+\delta) \Psi_T(\xi^{2(q-1)}) + C_{\delta}' \le (1+2\delta)\Psi_T(V(\xi)^2),
\end{multline*}
if $K$ is big enough. Then, the integral~\eqref{eq:integral} does not exceed
\[
(1+2\delta)\, \int_\R |\widehat{f}(\xi)|^2 \Psi_T(V(\xi)^2)\, {\rm d}\xi\,.
\]

Combining everything, we get, for big enough $K$,
\begin{multline*}
(1-\eps_1)\sigma^2 \int_\R |f(x)|^2 \Phi_T(U(x)^2)\, {\rm d}x \\
\le (1+2\delta)\, \int_\R |\widehat{f}(\xi)|^2 \Psi_T(V(\xi)^2)\, {\rm d}\xi
+ C_1(K) \cdot C_\eps  \theta \sigma^{2\theta-1}\,
\sum_{\la\in\Lambda'} |f(\la)|^2 U(\la)^{2\theta-1}\,.
\end{multline*}
Recalling that

\smallskip
$\eps_1$ and $\delta$ can be made arbitrarily small while $\sigma>1$ remains fixed,

\smallskip $U(\la)^{2\theta-1} = \max(K, |\la|)^{(2\theta - 1)(p-1)}$,

\smallskip
$\theta = sq = sp/(p-1)$, and therefore $(2\theta-1)(p-1) = 2sp -p+1 = (2s-1)p+1$,

\smallskip\noindent we get
\begin{multline*}
\int_\R |f(x)|^2 \Phi_T(U(x)^2)\, {\rm d}x
\le (1-\pho)\, \int_\R |\widehat{f}(\xi)|^2\Psi_T(V(\xi)^2)\, {\rm d}\xi \\
+ C\sum_{\lambda\in \Lambda'} |f(\lambda)|^2\max(K, |\lambda|)^{(2s-1)p + 1}
\end{multline*}
with some $\pho > 0$. Note that the integral in the RHS is finite since $f\in \Ho$ and the sum is finite by our assumption. Thus the LHS is finite as well.

Similarly, adding enough points to $M$ and repeating the same argument for the function $\widehat{f}$ with the function $\Psi_T$ instead of $\Phi_T$ we get
\begin{multline*}
\int_\R |\widehat{f}(\xi)|^2 \Psi_T(V(\xi)^2)\, {\rm d}\xi
\le (1-\pho)\, \int_\R |f(x)|^2\Phi_T(U(x)^2)\, {\rm d}x \\
+ C\sum_{\mu\in M'} |\widehat{f}(\mu)|^2\max(K, |\mu|)^{(2s-1)q + 1},
\end{multline*}
where the RHS is finite as we just established.

Summing these estimates and dividing by $\pho$ we get
\begin{multline*}
\int_\R |f(x)|^2\Phi_T(U(x)^2){\rm d}x + \int_\R |\widehat{f}(\xi)|^2\Psi_T(V(\xi)^2){\rm d}\xi \\
\le C\sum_{\lambda\in \Lambda'} |f(\lambda)|^2\max(K, |\lambda|)^{(2s-1)p + 1}
+ C\sum_{\mu\in M'} |\widehat{f}(\mu)|^2\max(K, |\mu|)^{(2s-1)q + 1}.
\end{multline*}

Finally, letting $T\to \infty$ and using the dominated convergence theorem we see that the LHS tends to
$$\int_\R |f(x)|^2 \max(K, |x|)^{2ps}{\rm d}x + \int_\R |\widehat{f}(\xi)|^2\max(K, |\xi|)^{2qs}\,
{\rm d}\xi,$$
which is at least $c||f||_{\Ho_{s, p, q}}^2$ for some $c > 0$. Bounding $\max(K, |y|)$ on the RHS
by $C'(1+|y|)$, we complete the proof of inequality \eqref{eq:*}.
\hfill $\Box$

\subsection{Removing finitely many points}
Our next goal is to show that the main estimate~\eqref{eq:*} holds for the pair
$(\Lambda, M)$.
For this,
we use the following lemma, cf. Duffin--Schaeffer~\cite[Lemma~IX]{Duffin-Schaeffer}:
\begin{lemma}[Duffin--Schaeffer]\label{hilbert lemma}
Let $\{v_n\}_{n = 0}^\infty$ be a system of vectors in a Hilbert space $H$.
Assume that there exists $C > 0$ such that, for all $v\in H$, we have
\begin{equation}\label{temp41}
||v||^2 \le C\sum_{n = 0}^\infty |\langle v, v_n\rangle |^2,
\end{equation}
and assume that the system $\{v_n\}_{n = 1}^\infty$ is complete in $H$.
Then there exists $C' > 0$ such that, for all $v\in H$, we have
\begin{equation}\label{temp411}
||v||^2 \le C'\sum_{n = 1}^\infty |\langle v, v_n\rangle |^2.
\end{equation}
\end{lemma}
For the reader's convenience, we provide a simple proof.

\medskip\noindent{\em Proof:}
We pick a small number $0 < \eps < 1$ to be determined later and take a finite linear combination of the vectors $v_1, v_2, \,  \ldots \,$, approximating $v_0$ with error at most $\eps$:
$$\bigl\| \sum_{n = 1}^M a_nv_n - v_0 \bigr\| \le \eps.$$
By Cauchy--Schwarz inequality we have
$$\bigl|\, |\langle v, \summ_{n = 1}^M a_n v_n\rangle| - |\langle v, v_0\rangle|\, \bigr| \le \eps \|v\|,$$
therefore,
\begin{equation}\label{temp42}
|\langle v, v_0\rangle|^2 - \bigl| \langle v, \summ_{n = 1}^M a_n v_n\rangle \bigr|^2
\le (2\| v_0 \| + \eps)\eps \| v \|^2 \le \eps(2\| v_0\|+1) \| v \|^2.
\end{equation}

We choose $\eps < (2C(2\| v_0\|+1))^{-1}$.
Juxtaposing estimates~\eqref{temp41}
and~\eqref{temp42}, we get
$$
\| v \|^2 \le 2C \bigl( \sum_{n = 1}^\infty |\langle v, v_n\rangle |^2 + |\langle v, \sum_{n = 1}^M a_n v_n\rangle|^2 \bigr)\,.$$
Once again applying the Cauchy--Schwarz inequality, we get
$$\| v\|^2 \le 2C\bigl( \sum_{n = 1}^\infty |\langle v, v_n\rangle |^2 + M\summ_{n = 1}^M |a_n|^2|\langle v, v_n\rangle|^2 \bigr),$$
which yields the estimate we need with $C' = 2C(1 + M\max |a_n|^2)$.
\hfill $\Box$

\medskip
By induction, the same conclusion holds  if we have an arbitrary finite set of vectors in place of $v_0$.

\medskip
With this lemma at hand, our next goal is to prove that, for $f\in\mathcal H$,
\begin{equation}\label{eq:temp*}
\| f \|_{\mathcal H_{s, p, q}}^2 \le C \,
\Bigl[\,
\sum_{\lambda\in\Lambda} (1+|\lambda|)^{(2s-1)p + 1} |f(\lambda)|^2
+ \sum_{\mu\in M} (1+|\mu|)^{(2s-1)q+1} |\widehat{f}(\mu)|^2\, \Bigr]\,.
\end{equation}
First, we show this assuming that $f\in\mathcal H_{s, p, q}$.
If $s\min(p, q)\ge 1$,  the point evaluations are continuous functionals on the Hilbert space
$\Ho_{s, p, q}$, therefore,
there exist vectors $\phi_x\in \Ho_{s, p, q}$ such that for all $f\in \Ho_{s, p, q}$ we have
$f(x) = \langle f, \phi_x\rangle$. Similarly, there exist vectors $\psi_\xi$ such that
$\widehat{f}(\xi) = \langle f, \psi_\xi\rangle$. Then, inequalty \eqref{eq:*} says that suitably normalized vectors $\phi_\lambda, \lambda\in\Lambda'$, and $\psi_\mu, \mu\in M'$, satisfy
assumption~\eqref{temp41} in Lemma~\ref{hilbert lemma}. On the other hand, since
$\Ho_{s, p, q}\subset \Ho$, by Theorem~1-UP, the system $\{\phi_\lambda\}_{\lambda\in \Lambda}\cup \{\psi_\mu\}_{\mu\in M}$ is complete in $\Ho_{s, p, q}$. Thus, by Lemma~\ref{hilbert lemma}, this system satisfies~\eqref{temp411}, which
gives us~\eqref{eq:temp*} for $f\in\mathcal H_{s, p, q}$.

Now, if $f\in\mathcal H$ and the sums over $\Lambda$ and $M$ in~\eqref{eq:temp*} are finite,
then the same sums over
$\Lambda'$ and $M'$ are finite as well (the sets $\Lambda'\setminus\Lambda$ and $M'\setminus M$ are finite, and, for $f\in\mathcal H$,
the point evaluations are bounded on the space and the Fourier sides).
Then, by Section~\ref{subsection:Adding}, $f\in\mathcal H_{s, p, q}$, and we arrive at the already considered case.

\subsection{Upper frame bound}
In this section we will show that, if the sequences $\Lambda$ and $M$ are $p$- and $q$-separated, respectively, then the converse to the estimate~\eqref{eq:temp*} is also true. Our first step is the following extension of Lemma~\ref{trace-claim}.
\begin{claim}\label{trace-claim-general}Let $\Gamma$ be a discrete set such that all gaps between consecutive elements of it have lengths at least $\delta$ for some $\delta > 0$. There exists an absolute constant $C > 0$ such that  for all $\theta\ge 1$ and for all $f\in H_\theta$ we have
$$
\delta\summ_{\gamma\in \Gamma} |f(\gamma)|^2 \le C \Bigl[\, \int_\R |f(x)|^2{\rm d}x
+ \delta^{2\theta} \int_\R |\xi|^{2\theta}|\widehat{f}(\xi)|^2\, {\rm d}\xi\Bigr].$$
\end{claim}

\noindent{\em Proof:}
If $\theta = 1$ then this follows from summing suitably scaled estimates from the Lemma~\ref{trace-claim} on all the intervals between consecutive elements of $\Gamma$.
For the general case, we have
\begin{multline*}
\int_\R |f(x)|^2{\rm d}x + \delta^2\int_\R |\xi|^2 |\widehat{f}(\xi)|^2
= \int_\R (1+ (\delta |\xi|)^2)|\widehat{f}(\xi)|^2\, {\rm d}\xi \\
\le 2 \int_\R (1+(\delta |\xi|)^{2\theta})|\widehat{f}(\xi)|^2\, {\rm d}\xi
= 2 \Bigl( \int_\R |f(x)|^2{\rm d}x
+ \delta^{2\theta} \int_\R |\xi|^{2\theta}|\widehat{f}(\xi)|^2\, {\rm d}\xi \Bigr).
\end{multline*}
Thus we need to multiply the constant $C$ by at most $2$.
\hfill $\Box$

\medskip
Now, we assume that the set $\Lambda$ is $p$-separated, that is, for all $x\in \R\backslash \Lambda$, the interval $I\subset \R\backslash \Lambda$ containing $x$ has length $\ge c(1+|x|)^{1-p}$.
As in the previous section, we fix $F\in C^\infty_0(-1, 1)$ with
$\displaystyle \int_\R |F(x)|^2\, {\rm d}x = 1$ and put $F_\upsilon(x) = F(x-\upsilon)$.
Applying Claim~\ref{trace-claim-general} with $\theta = sq$ and
$\delta = c(1+|\upsilon|)^{1-p}$ to the function $fF_\upsilon$, and multiplying both sides of the
resulting inequality by $(1+|\upsilon|)^{2sp}$, we get
\begin{multline}\label{upper-estimate}
(1+|\upsilon|)^{2sp+1-p}\sum_{\lambda\in\Lambda} |f(\lambda)|^2|F_\upsilon(\lambda)|^2\\
\le C\Bigl[ (1+|\upsilon|)^{2sp}\int_\R |f(x)|^2|F_\upsilon(x)|^2\, {\rm d}x +
\int_\R |\xi|^{2qs}|(\widehat{f}*\widehat{F}_\upsilon)(\xi)|^2\, {\rm d}\xi \Bigr],
\end{multline}
where we used the fact that for $x\in \supp F_\upsilon$, $(1+|x|)$ is comparable to
$(1+|\upsilon|)$. Now, we will integrate this inequality over $\upsilon\in\R$. For the LHS we get
$$
\int_\R (1+|\upsilon|)^{2sp+1-p}\sum_{\lambda\in\Lambda}
|f(\lambda)|^2|F_\upsilon(\lambda)|^2\, {\rm d}\upsilon\,,
$$
which is comparable to
$$
\sum_{\lambda\in \Lambda} (1+|\lambda|)^{(2s-1)p+1}|f(\lambda)|^2,
$$
since again, for $\lambda \in \supp F_\upsilon$, $(1+|\lambda|)$ is comparable to
$(1+|\upsilon|)$.

For the first term on the RHS of \eqref{upper-estimate}, by a similar computation, we get that its integral over $\upsilon\in \R$ is comparable to
$$\int_\R |f(x)|^2 (1+|x|)^{2sp}\, {\rm d}x.$$
Finally, for the second term on the RHS, its integral is
\begin{multline*}
\int_\R \int_\R |\xi|^{2qs}|(\widehat{f}*\widehat{F}_\upsilon)(\xi)|^2\, {\rm d}\xi {\rm d}\upsilon
= \int_\R\int_\R |\widehat{f}(\xi)|^2|\widehat{F}(\eta)|^2|\xi+\eta|^{2qs}\, {\rm d}\xi {\rm d}\eta \\
\le 2^{2qs} \int_\R\int_\R |\widehat{f}(\xi)|^2|\widehat{F}(\eta)|^2(|\xi|^{2qs}+|\eta|^{2qs})\,
{\rm d}\xi {\rm d}\eta
\le C \int_\R |\widehat{f}(\xi)|^2 (1+|\xi|^{2qs})\, {\rm d}\xi,
\end{multline*}
where in the first step, as before, we used the Plancherel theorem. Combining everything,
we see that the RHS of \eqref{upper-estimate} is bounded by $C||f||_{\Ho_{s, p, q}}^2$, that is,
$$\sum_{\lambda\in \Lambda} (1+|\lambda|)^{(2s-1)p+1}|f(\lambda)|^2 $$
is dominated by  $\|f\|_{\Ho_{s, p, q}}^2 $, as required.

Similarly, if $M$ is $q$-separated, we can get that
$$\sum_{\mu\in M}(1+|\mu|)^{(2s-1)q+1} |\widehat{f}(\mu)|^2 $$
is dominated by $ ||f||_{\Ho_{s, p, q}}^2 $.

Together, the results of this and two previous sections give us Theorem~\ref{thm-stable-uniq}.

\section{Non-uniqueness pairs}\label{sect:Non-uniq}
Here, we prove the second part of Theorem~\ref{thm-prelim}, which claims that any subcritical pair
$(\Lambda, M)$ is a non-uniqueness pair for the Schwarz space $\So$. As we have already mentioned,
we will prove slightly more:

\medskip\noindent{\bf Theorem~1-NUP.}\label{thm:1-NUP}
{\em Any subcritical pair $(\Lambda, M)$ is a non-uniqueness pair for the
Gelfand--Shilov space $\So(p, q)$, where $p$ and $q$ are the parameters from the
subcriticality assumption}.

It is worth mentioning that this theorem holds under much weaker assumptions on the pair $(\Lambda, M)$. We say that a sequence $\Gamma\subset\mathbb C$ has a density
with respect to the order $p$ if there exists a limit
\[
D(\Gamma, p) = \lim_{r\to\infty} \frac{n_\Gamma(r)}{r^p}\,.
\]
Here and elsewhere, $n_\Gamma(r) = |\Gamma \cap r\overline{\mathbb D}|$ denotes the counting function of the sequence $\Gamma$. In Theorem~1-NUP it suffices to assume that the sets $\Lambda_\pm = \Lambda \cap \mathbb R_\pm$, $M_\pm = M \cap \mathbb R_\pm$ have densities
\[
D(\Lambda_+, p), D(\Lambda_-,p) < \frac2p\,,
\]
\[
D(M_+, q), D(M_-, q) < \frac2q\,,
\]
and are $p$- and $q$-separated (in the sense of Definition~\ref{def:ser_seq}). Such a version
would require the full strength of Levin's theorem~\cite[Ch.~II, Theorem~5]{Levin}.

\subsection{Preliminaries}

\begin{definition}[{\sf $p$-smooth sequences}]
Let $p>0$ and let $\Gamma = (\gamma_j)$ be a sequence of points in $\mathbb C$ lying on a ray
${\rm arg\,}(\gamma_j)=\theta$, $|\gamma_1|<|\gamma_2|<\, \ldots $, $|\gamma_j|\uparrow\infty$.
We call the sequence $\Gamma$ $p$-smooth with a positive density $D=D(\Gamma, p)$ (with respect to the exponent $p$) if

\smallskip\noindent {\rm (a)}
\[
| n_\Gamma (r) - D r^p | = O(1), \quad r\to\infty\,;
\]

\smallskip\noindent {\rm (b)} there exists $d>0$ such that
$|\gamma_{j+1}|-|\gamma_j| \ge d(1+|\gamma_j|)^{1-p}$, $j=1, 2, \ldots $.
\end{definition}

\begin{claim}\label{lemma-enlarging}
Let $p>0$, and let $\Gamma = (\gamma_j) \subset \R_+$, $\gamma_1 < \gamma_2 < \ldots$,
$\gamma_j\uparrow\infty$.
Suppose that
\begin{equation}\label{eq:***}
\liminf_{j\to\infty} \gamma_j^{p-1} (\gamma_{j+1}-\gamma_j) > \delta >0\,.
\end{equation}
Then, for any $D > 1/(p\delta)$, there exists a $p$-smooth supsequence
$\Gamma' \supset \Gamma$ with density~$D$, and $|\Gamma'\setminus \Gamma|=\infty$.
\end{claim}

\noindent{\em Proof of Claim~\ref{lemma-enlarging}}:
First, we note that it suffices to prove the claim for $p=1$. The general case follows by applying the case $p=1$ to the sequence $(\gamma_j^{p})$.

Split $[0, +\infty)$ into intervals of length $1/D<\delta$. By our assumption,
each of these intervals, except maybe finitely many ones, contains at most one point
of the set $\Gamma$. At the same time, there are infinitely many intervals which do not contain points of $\Gamma$ (otherwise, $\liminf r^{-1} n_\Gamma (r) \ge D$, while \eqref{eq:***} yields that $\limsup r^{-1}n_\Gamma (r) \le \delta^{-1}$).
If the interval does not contain a point of $\Gamma$, we place a new point at its center. \mbox{}\hfill $\Box$

\medskip
In what follows, relying on Claim~\ref{lemma-enlarging},
we assume that $p, q > 1$, $\displaystyle \tfrac1p + \tfrac1q = 1$,
and add to each of the sets $\Lambda_\pm = \Lambda \cap \mathbb R_\pm$, $M_\pm = M \cap \mathbb R_\pm$ infinitely many points so that the enlarged sets $\Lambda'_\pm$ and $M'_\pm$ become
$p$- and $q$-smooth and satisfy
\[
D(\Lambda'_+, p) = D(\Lambda'_-,p) < \frac2p\,,
\]
\[
D(M'_+, q) = D(M'_-, q) < \frac2q\,.
\]

\subsection{The main lemma}
The proof of Theorem~1-NUP relies on the construction of the sequences of entire
functions $(\Phi_\lambda)_{\lambda\in\Lambda'}$ and $(\Psi_\mu)_{\mu\in M'}$ that interpolate the
$\delta$-functions $\delta_\lambda$ and $\delta_\mu$ and satisfy certain bounds uniformly with respect to $\lambda$ and $\mu$.

\begin{lemma}\label{lemmaA}
Let $p > 1$, $q=p/(p-1)$. Let $\Lambda'\subset\mathbb R$, and let
the sets $\Lambda'_\pm$ be $p$-smooth with
\[
D(\Lambda'_+, p) =  D(\Lambda'_-, p) < \frac2p\,.
\]
Then, for any sufficiently small $a>0$, there exist $C<\infty$, $a'>a>a''>0$, and a sequence
of entire functions $(\Phi_\lambda)_{\lambda\in\Lambda'}$  with the following properties
\begin{align}
\label{Phi_1}
\Phi_\lambda(\lambda') &=
\begin{cases}
1, & \lambda'=\lambda,  \\
0, & \lambda'\in\Lambda'\setminus \{\lambda\};
\end{cases}
\\
\label{Phi_2}
|\Phi_\lambda(x)| &\le Ce^{-a''|x|^p + a|\lambda|^p}, \quad x\in\mathbb R;
\\
\label{Phi_3}
|\widehat{\Phi}_\lambda(\xi)| &\le C e^{-a'|\xi|^q+a|\lambda|^p},
\quad \xi\in\mathbb R\,.
\end{align}
\end{lemma}

Note that simultaneously with the sequence  $(\Phi_\lambda)_{\lambda\in\Lambda'}$ the lemma also provides us with the
sequence  of entire functions $(\Psi_\mu)_{\mu\in M'}$ such that
\begin{align}
\label{Psi_1}
\Psi_\mu(\mu') &=
\begin{cases}
1, & \mu'=\mu,  \\
0, & \mu'\in M'\setminus \{\mu\};
\end{cases}
\\
\label{Psi_2}
|\Psi_\mu(\xi)| &\le Ce^{-a''|\xi|^q + a|\mu|^q}, \quad \xi\in\mathbb R;
\\
\label{Psi_3}
|\widecheck{\Psi}_\mu(x)| &\le C e^{-a'|x|^p+a|\mu|^q},
\quad x\in\mathbb R,
\end{align}
where $\widecheck{\Psi}$ denotes the inverse Fourier transform.
We use bounds~\eqref{Phi_2} and~\eqref{Psi_2} for the functions $\Phi$ and $\Psi$ with the
same positive constant $a''$ (otherwise, we take the minimum of two constants).
Similarly, we assume that the positive constant $a'$ in~\eqref{Phi_3} and in~\eqref{Psi_3}
is the same.

\medskip
The proof of Lemma~\ref{lemmaA} relies on several notions and results pertaining to the classical theory of entire functions, in particular, on the B.~Ya.~Levin estimates for entire functions with regular zeroes.

We proceed with the proof of Theorem~1-NUP, assuming the existence of the sequences  $(\Phi_\lambda)$ and $(\Psi_\mu)$, and then will provide the proof of Lemma~\ref{lemmaA}.

\subsection{Free interpolation by Fourier pairs}\label{subsect-FreeInterpol}
We fix $\delta\le 1/(4C)$, where $C$ is a constant from estimates \eqref{Phi_2}, \eqref{Phi_3}, \eqref{Psi_2}, and \eqref{Psi_3}, and choose $L$ so large that
\[
\sum_{\lambda\in\Lambda'\setminus [-L, L]} e^{-(a'-a)|\lambda|^p} + \sum_{\mu\in M'\setminus [-L, L]} e^{-(a'-a)|\mu|^q} < \delta.
\]
Then, we set $\Lambda_L'=\Lambda'\setminus [-L, L]$, and $M_L'=M'\setminus [-L, L]$. By $\sf B$ we denote the Banach space of pairs of fast decaying sequences
\[
\kappa = (\alpha, \beta), \quad  \alpha\colon \Lambda_L'\to\mathbb C, \quad  \beta\colon M_L'\to\mathbb C,
\]
endowed with the weighted $\ell^1$ norm
\[
\| \kappa \| \stackrel{\rm def}= \sum_{\lambda\in\Lambda_L'} |\alpha(\lambda)|e^{a|\lambda|^p} + \sum_{\mu\in M_L'} |\beta(\mu)| e^{a|\mu|^q}\,.
\]

\begin{lemma}\label{lemma-Claim}
For any sequence $\kappa\in B$, there exists an entire function $f$ from the Gelfand--Shilov space
$\So(p, q)$ that solves the interpolation problem
\begin{equation}\label{eqI.P.}
\begin{cases}
f(\lambda) = \alpha(\lambda), \quad &\lambda\in\Lambda_L', \\
\widehat{f}(\mu) = \beta (\mu), \quad & \mu\in M_L'\,.
\end{cases}
\end{equation}
\end{lemma}

Note that there is no hope to solve such an interpolation problem (even with the Schwartz function $f$)
with the full set of nodes $\Lambda' \cup M'$. Evident obstacles are
classical, as well as non-classical, Poisson summation formulas. It is curious that removing only finitely many nodes from both sets, we get rid of the hurdles.
This remark can be rephrased in terms of crystalline measures. Lemma~\ref{lemma-Claim} yields that the only locally finite measure of tempered growth supported on $\Lambda'$ such that its Fourier transform is a locally finite measure of tempered growth supported on $M'$ is a zero measure. Since $\Lambda\backslash \Lambda'$ and $M\backslash M'$ are finite,  by a simple linear algebra argument we can conclude that the vector space of such measures, but supported on $\Lambda$ and $M$ respectively, is finite-dimensional.

\subsection{Proof of Lemma~\ref{lemma-Claim}}
We construct the interpolating entire function $f$ by an iterative procedure. We start with
\[
f_1 = \sum_{\lambda'\in \Lambda_L'} \alpha(\lambda') \Phi_{\lambda'} + \sum_{\mu'\in M_L'} \beta(\mu') \widecheck{\Psi}_{\mu'},
\]
and estimate the sizes of $f_1$ and $\widehat{f_1}$. We have
\begin{align*}
|f_1(x)| &\le \sum_{\lambda'\in \Lambda_L'} |\alpha(\lambda')|\, |\Phi_{\lambda'}(x)|
+ \sum_{\mu'\in M_L'} |\beta(\mu')|\, |\widecheck{\Psi}_{\mu'}(x)| \\
&\le C\,\Bigl[\,
\sum_{\lambda'\in \Lambda_L'} |\alpha(\lambda')|\, e^{-a''|x|^p + a|\lambda'|^p}
+ \sum_{\mu'\in M_L'} |\beta(\mu')|\, e^{-a'|x|^p + a|\mu'|^q}
\Bigr] \quad ({\text{by}\ \eqref{Phi_2}, \eqref{Psi_2}})
\\
&\le C \| \kappa \| e^{-a''|x|^p} \qquad (\text{since}\ a''<a'),
\end{align*}
and, similarly,
\[
|\widehat{f_1}(\xi)| \le C\| \kappa\| e^{-a''|\xi|^q}\,.
\]
That is, the function $f_1$ lies in the Gelfand--Shilov space $\mathcal S(p, q)$.

Next, we look at the differences
\[
\alpha_1(\lambda') = f_1(\lambda') - \alpha(\lambda') = \sum_{\mu'\in M'_L} \beta(\mu') \widecheck{\Psi}_{\mu'}(\lambda'),
\qquad \lambda'\in \Lambda_L',
\]
and
\[
\beta_1(\mu') = \widehat{f_1}(\mu') - \beta(\mu') = \sum_{\lambda'\in\Lambda_L'} \alpha(\lambda') \widehat{\Phi}_{\lambda'}(\mu'),
\qquad \mu'\in M_L'.
\]
We have
\begin{align*}
\sum_{\lambda'\in\Lambda_L'} |\alpha_1(\lambda')|\, e^{a|\lambda'|^p}
&\underset{\eqref{Psi_3}}\le C\, \sum_{\mu'\in M_L'} |\beta(\mu')| e^{a|\mu'|^q}\, \cdot \,
\sum_{\lambda'\in\Lambda_L'} e^{(a-a')|\lambda'|^p}
\\
&\le
C\delta\, \sum_{\mu'\in M_L'} |\beta(\mu')| e^{a|\mu'|^q}  \qquad (\text{choice\ of\ } L)
\\
&\le C\delta\|\kappa\| \\
&\le\tfrac14\,\|\kappa\| \qquad (C\delta\le\tfrac14).
\end{align*}
Similarly,
\[
\sum_{\mu'\in M_L'} |\beta_1(\mu')|\, e^{a|\mu'|^q} \le \tfrac14\, \| \kappa \|\,.
\]

Now, we let $\kappa_1 = (\alpha_1, \beta_1)$ and repeat the construction, noting that $\| \kappa_1 \| \le \tfrac12 \| \kappa \|$.
Clearly, the process will converge with the speed of a geometric progression to the solution of the interpolation problem~\eqref{eqI.P.}, which also belongs to the space $\mathcal S(p, q)$.
This proves Lemma~\ref{lemma-Claim}. \hfill $\Box$

\subsection{Completing the proof of Theorem~1-NUP}
To finish the proof of the theorem, we consider sequences
$\kappa=(\alpha, \beta)$ with $\alpha$ vanishing on
$\Lambda\setminus [-L, L]$ and $\beta$ vanishing on $M\setminus [-L, L]$.
Since the set $\Lambda'\backslash (\Lambda\bigcup [-L, L])$ is infinite,
Lemma~\ref{lemma-Claim} will provide us
with an infinite-dimensional subspace of functions $f\in\So(p, q)$ such that
$f\big|_{\Lambda\setminus [-L, L]} = 0$ and
$\widehat{f}\,\big|_{M\setminus [-L, L]} = 0$. By a simple linear algebra argument, this subspace contains non-zero  functions $f$ satisfying finitely many
linear equations
\[
\begin{cases}
f(\lambda)=0, \quad &\lambda\in\Lambda\cap [-L, L], \\
\widehat{f}(\mu)=0, \quad &\mu\in M\cap [-L, L].
\end{cases}
\]
This completes the proof of Theorem~1-NUP modulo the proof of Lemma~\ref{lemmaA}. \hfill $\Box$

\subsection{Proof of Lemma~\ref{lemmaA}}
We will be using  several notions and results from the classical theory of entire functions, see Levin~\cite{Levin, Levin2}.

\subsubsection{The class $\mathsf K_p$}

\begin{definition}[The class $\mathsf K_p$]
A $2\pi$-periodic continuous function $k$ belongs to the class
$\mathsf K_p$ ($p$ is a positive constant) if
\[
k'' + p^2 k = \sum_{j=1}^n m_j \delta_{\theta_j}\,,
\]
where $\delta_\theta$ is a unit point measure at $\theta$, the masses $m_j$ are positive, and
the differential operator on the LHS is understood distributionally.
\end{definition}
Formally the latter means that, for any $2\pi$-periodic $C^\infty$-function $\phi$,
\begin{equation}\label{eq:K_p}
\int_{\mathbb R /2\pi\mathbb Z} k(\phi''+p^2 \phi) = \sum_{j=1}^n m_j \phi(\theta_j).
\end{equation}

We assume that the points $\theta_1, \theta_2, \, \ldots \,, \theta_n \in \R/2\pi\mathbb Z$
are arranged counterclockwise
and mention two simple facts concerning functions $k\in \mathsf K_p$.
\begin{itemize}
\item[(i)] On each interval $(\theta_j, \theta_{j+1})$, the function $k$ is $p$-trigonometric, i.e., $k(\theta) = a_j \cos p\theta + b_j\sin p\theta$.
\item[(ii)] The masses $m_j$ are equal to the jumps of the derivative of $k$, i.e., $m_j = k'(\theta_j+0)-k'(\theta_j-0)$.
\end{itemize}

For the reader's convenience, we sketch the proofs of these facts\footnote{
Functions from $\mathsf K_p$ are $p$-trigonometrically convex~\cite[Ch.I, \S16]{Levin}, \cite[Lecture~8]{Levin2}, and these facts are special cases of the corresponding properties of $p$-trigonometrically convex functions.}.

\smallskip\noindent (i) If the function $k$ is $C^\infty$-smooth on $(\theta_j, \theta_{j+1})$,
then, taking an arbitrary $\phi\in C_0^\infty(\theta_j, \theta_{j+1})$ and integrating by parts
twice on the LHS of~\eqref{eq:K_p}, we conclude that $k'' + p^2 k = 0$
on $(\theta_j, \theta_{j+1})$ and, therefore, $k$ is $p$-trigonometric therein.
If $k$ is only continuous, then we take a mollifier $\psi\in C_0^\infty (-\pi, \pi)$ with
integral $1$, put $\psi_t = t \psi (t\, \cdot\,)$, apply the previous argument to the convolution $k \ast \psi_t$, conclude that $k\ast\psi_t$ is $p$-trigonometric on $\bigl( \theta_j + \pi t^{-1}, \theta_{j+1}- \pi t^{-1} \bigr)$ for each $t$, and then let $t\to\infty$, noting that if a sequence
of $p$-trigonometric functions converges pointwise, then the limiting function is $p$-trigonometric
as well.

\smallskip\noindent (ii) The integral on the LHS of~\eqref{eq:K_p} equals
\[
\sum_{j=1}^n \int_{\theta_j}^{\theta_{j+1}} k(\phi'' + p^2 \phi)
= \sum_{j=1}^n (k'(\theta_j+0) - k'(\theta_j-0)) \phi (\theta_j)
\]
(we integrated twice by parts and combined the terms that appeared with opposite signs).

\subsubsection{Entire functions with smooth sequences of zeroes}
With each function $k\in {\sf K}_p$ we associate a class of $k$-smooth discrete sets in the plane.
\begin{definition}[$k$-smooth sets]
A discrete set of points
\[
Z = \bigcup_{j=1}^n Z_j, \quad Z_j = \bigl\{
z_{\ell, j} = |z_{\ell, j}|e^{{\rm i}\theta_j}
\colon \ell \ge 1 \bigr\}\,,
\]
is {\em $k$-smooth}, if it satisfies the following conditions:

\smallskip\noindent{\rm (a)} for each $1\le j \le n$,
the set $Z_j$ is $p$-smooth on the ray $\{{\rm arg}(z) = \theta_j\}$ with density $D_j=(2\pi p)^{-1}m_j$,
that is,
\[
\bigl| Z_j \cap [0, re^{{\rm i}\theta_j} ] \bigr| = \frac{m_j}{2\pi p}\, r^p + O(1),
\quad r\to\infty;
\]

\smallskip\noindent{\rm (b)}
the disks $D_{\ell, j}= D\bigl( z_{\ell, j}, d(1+|z_{\ell, j}|)^{1-p} \bigr)$,
$\ell\ge 1$, $1\le j \le n$, are disjoint for some $d>0$.
\end{definition}

\smallskip
The following theorem is a very special case of Levin's theorem~\cite[Ch.~II, Theorem~5]{Levin}.

\noindent{\bf Theorem~L}
{\em Let $k\in \mathsf K_p$. Then, for any $k$-smooth set $Z$,
there exists an entire function $S$, whose zeroes are simple and coincide with $Z$, such that,
for every $\varepsilon>0$,}
\begin{align*}
\tag{I}
|S(re^{{\rm i}\theta})| &\le C_\varepsilon e^{(k(\theta)+\varepsilon)r^p},
\qquad {\rm everywhere\ in\ } \mathbb C,
\\
\tag{II}
|S(re^{{\rm i}\theta})| &\ge c_\varepsilon  e^{(k(\theta)-\varepsilon)r^p},
\qquad {\rm everywhere\ in\ } \mathbb C\setminus \bigcup_{j=1}^n\,\bigcup_{\ell\ge 1} D_{\ell, j}.
\end{align*}

Note that estimate (II) yields the lower bound
\[
\tag{III}
|S'(z_{\ell, j})| \ge c_\varepsilon e^{(k(\theta_j)-2\varepsilon)|z_{\ell, j}|^p}\,.
\]
To get it, one has to apply the minimum modulus principle to the harmonic function
$\log |S(z)/(z-z_{\ell, j})|$ in the disk $D_{\ell, j}$.

\medskip
For the reader's convenience, we will provide in Appendix a short self-contained proof of Theorem~L,
which does not use the machinery developed in~\cite[Ch.~II]{Levin}.

\subsubsection{Choice of the functions $k_p(\theta)$}
We will apply Theorem~L to a special family of functions $k_p\in {\sf K}_p$.
These functions will possess the following properties:

\smallskip\par\noindent (a)
they are symmetric with respect to $\theta=0$ and $\displaystyle \theta=\pm\frac{\pi}2$, that is,
$k_p(-\theta)=k_p(\theta)$, and $k_p(\pi-|\theta|)=k_p(\theta)$ for $-\pi \le\theta \le \pi$;

\smallskip\par\noindent (b)
$k_p(0) = k_p(\pm\pi)<0$;

\smallskip and

\smallskip\par\noindent (c)
the functions $k_p$ are increasing and $p$-trigonometric on $\displaystyle \bigl[ 0, \frac{\pi}{2p}\, \bigr]$.

\medskip
For $1<p\le 2$, we set
\begin{equation}\label{relation}
k_p(\theta) = \alpha \sin p\theta - \beta \cos p\theta, \quad 0 \le \theta \le \frac{\pi}2,
\end{equation}
with $\alpha, \beta >0$, and then extend it to the whole interval $[-\pi, \pi]$ using the symmetries
and letting
\[
k_p(\pi-\theta)=k_p(-\pi+\theta)= k_p(-\theta)=k_p(\theta), \qquad 0\le \theta \le \frac{\pi}2\,.
\]
To keep positivity of the jumps of $k'$ at $\theta = \pm \pi/2$, we require that $k_p'(\pi/2-0)\!<\!0$,
that is, for $1<p<2$, $\beta < \alpha |\cot(\pi p/2)|$ (for $p=2$, this requirement
is void).

For $p>2$, we let $p=2^np'$ with $n\in\mathbb N$ and $1<p'\le 2$, and define
$k_p(\theta)=k_{p'}(2^n\theta)$. Note that relation~\eqref{relation}
continues to hold for $0\le \theta \le \pi/(2p)$.

\medskip
To adjust the function $k_p$ to the densities $D(\Lambda'_\pm, p)$, we choose $\sigma <1$ so that
\[
D(\Lambda'_+, p) = D(\Lambda'_-, p) = \frac{2\sigma}p\,,
\]
and let $\alpha = 2\pi \sigma /p$. Then,
\begin{align*}
k_p'(+0)-k_p'(-0) &= k_p'(\pm \pi +0) - k'_p(\pm\pi -0) \\
&= 2\alpha p =  2\pi p \cdot \frac{2\sigma}p = 2\pi p D(\Lambda'_\pm, p),
\end{align*}
which will be needed momentarily for the application of Levin's theorem.
The small parameter $\beta$ in the definition of the function $k_p$ will be fixed a bit later.

\subsubsection{The entire function $\Phi$}
We apply Theorem~L with the function $k_p$ and a set $Z\supset\Lambda'$.

We enlarge $\Lambda'$, adding to it $p$-smooth sets lying on the rays corresponding to discontinuities of $k_p'$ different from $\theta=0$ and $\theta=\pm\pi$, and
choosing their densities according to the jumps
of $k_p'$ at these discontinuities, so that the new set $Z$ meets
conditions (a) and (b) from the definition of the $k$-smooth set.
Then Theorem~L provides us with an entire function $\Phi$ with simple zeroes at $Z$, such that
for every $\varepsilon>0$,
\[
|\Phi(re^{{\rm i}\theta})| \le C_\varepsilon e^{(k(\theta)+\varepsilon)r^p},
\qquad {\rm everywhere\ in\ } \mathbb C,
\]
and
\[
|\Phi'(\lambda)| \ge c_\varepsilon  e^{-(\beta+\varepsilon)|\lambda|^p},
\qquad \lambda\in\Lambda'\,,
\]
where $\beta=-k_p(0)=-k_p(\pm \pi)$.

\subsubsection{Decay of the Fourier transform}
The next step is to estimate the decay of the Fourier transform of the function $F_\lambda (z) = \Phi(z)/(z-\lambda)$, uniformly in $\lambda\in\Lambda'$. First, we claim that, for every $\varepsilon>0$,
\begin{equation}\label{eqF}
|F_\lambda (re^{{\rm i}\theta})| \le C_\varepsilon e^{(k_p(\theta)+\varepsilon)r^p}\,,
\qquad {\rm everywhere\ in\ } \mathbb C\,.
\end{equation}
Clearly, this bound holds in $\{z\colon |\Im (z)|\ge 1\}$. Then, by the Phragm\'en--Lindel\"of theorem, it also holds in the strip $\{z\colon |\Im (z)|\le 1\}$ (recall that the entire function $F_\lambda$ has finite order of growth).

Shifting the integration line in the definition of the Fourier transform $\widehat{F}_\lambda$, we get
\[
\widehat{F}_\lambda(\xi) =
\int_\mathbb R F_\lambda(x+{\rm i}y) e^{-2\pi {\rm i} \xi (x+{\rm i}y)}\, {\rm d} x\,,
\]
with any $y\in\mathbb R$. The shift is justified since the function $k_p(\theta)$ stays negative in the angles $\{|\theta|\le\delta\}$ and $\{|\pi - \theta|\le\delta\}$ with some $\delta>0$, and therefore,
$|F_\lambda(z)|$ decays in these angles like $e^{-c|z|^p}$. Thus, for every $\varepsilon>0$,
\begin{equation}\label{Fourier contour shift}
\bigl| \widehat{F}_\lambda (\xi) \bigr| \le
C_\varepsilon\, \inf_{y\in\mathbb R}\, \sup_{x\in \mathbb R}
|F_\lambda(x+{\rm i}y)| e^{\varepsilon |x|^p + 2\pi \xi \cdot y}\,.
\end{equation}
We will choose $y$ satisfying $\xi \cdot y \le 0$, so that the last factor in \eqref{Fourier contour shift} is always $e^{-2\pi |\xi|\cdot |y|}$.

Letting $x+{\rm i}y = re^{{\rm i}\theta}$ and taking into account estimate~\eqref{eqF},
in order to proceed,
we need to bound the expression $(k_p(\theta)+2\varepsilon)r^p$ for $-\pi \le\theta\le\pi$. Because of the symmetries of $k_p$, we can restrict ourselves to the region $0\le \theta \le \pi/2$.

\begin{claim}\label{claim:trig-ineq}
For all $s\in (\sigma^q, 1)$, there exists $b_0(s)$ such that, for all $\theta\in \displaystyle \bigl[ 0, \frac{\pi}2\, \bigr]$
and all $b>b_0(s)$, we have
\begin{equation}\label{eq:ineq**}
k_p(\theta) < 2\pi \frac{b^p}p\, \sin^p\theta,
\end{equation}
provided that
\[
\beta = 2\pi s \frac{b^{-q}}q\,.
\]
\end{claim}

\noindent{\em Proof:}
Choose $\theta_0\in \bigl( 0, (2p)^{-1}\pi \bigr) $ so that
$s \cos p\theta_0 = \sigma^q$. Note that, for $\theta_0\le \theta \le \pi/2$,
the inequality holds since, increasing $b$, we unboundedly increase the minimum of the RHS of~\eqref{eq:ineq**}, while the LHS of~\eqref{eq:ineq**} stays bounded.

For $0\le \theta \le \theta_0$, we have
$ k_p(\theta) = \displaystyle \frac{2\pi\sigma}p \sin\theta - \beta\cos\theta $,
so we need to show that
\[
\frac{2\pi \sigma}p\, \sin p\theta < \beta \cos p\theta + 2\pi \frac{b^p}p\, \sin^p\theta\,,
\]
or (substituting the value of $\beta$) that
\[
\frac{\sigma}p\, \sin p\theta < s \frac{b^{-q}}q\, \cos p \theta
+ \frac{b^p}p\, \sin^p\theta\,.
\]
Then, since $\cos p\theta \ge \cos p\theta_0$, it suffices to show that
\[
\frac{\sigma}p\, \sin p\theta < \sigma^q\, \frac{b^{-q}}q\, + \frac{b^p}p\, \sin^p\theta\,.
\]
By the Young inequality, the RHS is not less than $\sigma\sin\theta$, while, by the concavity of the sine function on
$\displaystyle \bigl[0, \frac{\pi}2\, \bigr]$, for $p>1$, $\sin p\theta < p\sin\theta$ everywhere on
$\displaystyle \bigl( 0, \frac{\pi}{2p}\, \bigr]$. This proves
estimate~\eqref{eq:ineq**} in the range $0 < \theta \le \theta_0$. For $\theta=0$, the estimate holds
because $k_p(0)<0$. \hfill $\Box$

\medskip
From now on, we fix $s\in (\sigma^q, 1)$, and let $\beta = 2\pi s b^{-q}/q$, with some $b>b_0(s)$ to be specified later. Then, by the claim,
\[
k_p(\theta) + 2\varepsilon \le 2\pi \frac{b^p}p\, \sin^p\theta.
\]
Note that $\eps > 0$ can be chosen independent of $\theta$ by continuity of the functions $k_p$ and sine. Hence,
\[
(k_p(\theta) + 2\varepsilon)r^p - 2\pi |\xi| \cdot |y| \le 2\pi\Bigl[\,
\frac{b^p}p\, |y|^p - |\xi| \cdot |y| \,\Bigr]\,.
\]
Thus, for $x+{\rm i}y=re^{{\rm i}\theta}$,
\begin{multline*}
|F_\lambda(x+{\rm i}y)| e^{\eps |x|^p-2\pi |\xi|\, |y|}
\le C_\eps e^{(k_p(\theta)+\eps)|x+{\rm i}y|^p +\eps |x|^p - 2\pi |\xi|\, |y|} \\
\le C_\eps e^{(k_p(\theta)+2\eps)|x+{\rm i}y|^p  - 2\pi |\xi|\, |y|}
\le C_\eps e^{2\pi [b^p|y|^p/p - |\xi|\, |y|]}\,.
\end{multline*}
Minimizing the exponent on the RHS, we choose
\[
|y| = \Bigl( \frac{|\xi|}{b^p} \Bigr)^{1/(p-1)}\,.
\]
Then,
\[
2\pi \Bigl[\,
\frac{b^p}p\, |y|^p - |\xi| \cdot |y| \,\Bigr] = -2\pi \frac{b^{-q}}q\, |\xi|^q\,,
\]
and therefore,
\[
\bigl| \widehat{F}_\lambda (\xi) \bigr| \le C_\varepsilon' e^{-\beta' |\xi|^q}\,,
\qquad \xi\in\mathbb R, \ \lambda\in\Lambda',
\]
with
\[
\beta' = 2\pi \frac{b^{-q}}q = \frac{\beta}s\,.
\]

\subsubsection{Completing the proof of Lemma~\ref{lemmaA}}

It remains to let
\[
\Phi_\lambda (z) = \frac{\Phi(z)}{\Phi'(\lambda)(z-\lambda)} = \frac{F_\lambda(z)}{\Phi'(\lambda)}\,,
\qquad \lambda\in\Lambda'\,.
\]
Evidently,
\[
\Phi_\lambda (\lambda') =
\begin{cases}
1, &\lambda'=\lambda, \\
0, &\lambda'\in\Lambda'\setminus\{\lambda\}\,.
\end{cases}
\]
Furthermore, for every $\varepsilon>0$, we have
\[
|\Phi_\lambda (x) |
\le C_\varepsilon e^{-(\beta-\varepsilon)|x|^p+(\beta+\varepsilon)|\lambda|^p}\,
\qquad x\in\mathbb R, \ \lambda\in\Lambda',
\]
and
\[
|\widehat{\Phi}_\lambda(\xi)| \le C_\eps e^{-\frac{\beta}s\,|\xi|^q + (\beta+\varepsilon)|\lambda|^p}\,,
\qquad \xi\in\mathbb R, \ \lambda\in\Lambda'\,.
\]

At last, given sufficiently small $a>0$, we choose $\beta$ so that $\beta<a<\beta/s$
(recall that we can choose $\beta$ as small, as we please, choosing the parameter $b$ in
Claim~\ref{claim:trig-ineq} sufficiently large), and then choose $\varepsilon< \min(a-\beta, \beta)$. Letting $a'=\beta/s$ and $a''=\beta-\varepsilon$, we obtain the needed
upper bounds:
\[
|\Phi_\lambda (x) |
\le C e^{-(\beta-\varepsilon)|x|^p+(\beta+\varepsilon)|\lambda|^p}
= C e^{-a''|x|^p+a|\lambda|^p}\,, \quad x\in\mathbb R, \ \lambda\in\Lambda',
\]
and
\[
|\widehat{\Phi}_\lambda(\xi)| \le Ce^{-\frac{\beta}s\,|\xi|^q + (\beta+\varepsilon)|\lambda|^p}
= Ce^{-a'|\xi|^q + a|\lambda|^p}\,,
\quad \xi\in\mathbb R, \ \lambda\in\Lambda'\,,
\]
thus completing the proof of Lemma~\ref{lemmaA}. \hfill $\Box$

\section*{Acknowledgement}
The work of Aleksei Kulikov was supported by Grant 275113 of the Research Council of Norway,
BSF Grant 2020019, ISF Grant 1288/21, and by The Raymond and Beverly Sackler Post-Doctoral Scholarship.
The work of Fedor Nazarov was supported in part by the NSF grant DMS-1900008 and BSF grant 2020019.
The work of Mikhail Sodin was supported by ERC Advanced Grant 692616, BSF Grant 2020019, and by ISF Grant 1288/21.

\section*{Appendix: Proof of Theorem~L}
\renewcommand{\theequation}{A.\arabic{equation}}
\setcounter{equation}{0}
\renewcommand{\theclaim}{A.\arabic{claim}}
\setcounter{claim}{0}
\renewcommand{\thelemma}{A.\arabic{lemma}}
\setcounter{lemma}{0}

\subsection*{Brelot--Hadamard representation of the function $K$}

Consider the function $K(re^{{\rm i}\theta}) = k(\theta)r^p$, where
$k$ is a function from the class $\mathsf K_p$, which is $p$-trigonometric
outside the points $\{\theta_j\}$, $1\le j \le n$. As before, we assume that these
points are arranged counterclockwise.
The function $K$ is a continuous subharmonic function
in $\mathbb C$ with the Riesz measure
\[
\mu =
\frac1{2\pi}\, \Delta K = \frac1{2\pi p}\, \sum_{j=1}^n m_j \delta_{\theta_j} \otimes {\rm d}(r^p)\,,
\]
where $m_j = k'(\theta_j+0)-k'(\theta_j-0) > 0$ are the jumps of $k'$.
That is, for any compactly supported smooth test function $f$, we
have\footnote{To verify this, one can apply Green's formula to $K$ and $f$ in the truncated sectors
$\Omega_{\rho, R, j} = \{re^{{\rm i}\theta}\colon
\rho < r < R, \theta_j < \theta < \theta_{j+1} \}$ with sufficiently large $R$, add the results,
and then let $\rho\to 0$.}
\[
\int_{\mathbb C} K \Delta f =
\sum_{j=1}^n m_j \int_0^\infty f(re^{{\rm i}\theta}) r^{p-1}\, {\rm d}r\,.
\]
In particular, the function $K$ is harmonic in each angle $\{\theta_j < \arg (z) < \theta_{j+1}\}$.
In addition, if $p$ is a positive integer, then
\begin{equation}\label{eq:Lindelof}
\textstyle \sum_j m_j e^{{-{\rm i}p\theta_j}} = 0.
\end{equation}
To see this, one can apply integration by parts in each interval
$(\theta_j, \theta_{j+1})$ to the functions $k(\theta)$ and $e^{-{\rm i}p\theta}$ and add the results up.

\medskip
Set $\mu_R = \done_{R\overline{\mathbb D}}\, \mu$. Let
\[
\mathsf E_s (w) = (1-w)e^{w+w^2/2+\, \ldots\, + w^s/s}, \quad s\in \mathbb N\,,
\]
be the Weierstrass factor, and let
$ \mathsf H_s(w) = \log |\mathsf E_s(w)|$.
Choose $s < p \le s+1$ and set
\[
U_R(z) = \int \mathsf H_s (z/\zeta)\, {\rm d}\mu_R (\zeta)\,.
\]
This is a subharmonic function with $\Delta U_R = (2\pi)^{-1} \mu_R$
(where the Laplacian is understood in the sense of distributions).

\begin{lemma}[A special case of the Brelot--Hadamard theorem] \mbox{}

\smallskip\noindent
1. We have $ U_R(z) \le C|z|^p$, $z\in\mathbb C$,
uniformly in $R\ge 1$.

\smallskip\noindent
2. The family $U_R$ converges as $R\to\infty$, and
\[
\lim_{R\to\infty} U_R(z) =
\begin{cases}
K(z), & p\notin \mathbb N, \\
K(z) + {\rm Re}(cz^p), & p\in\mathbb N,
\end{cases}
\]
locally uniformly in $\mathbb C$.
\end{lemma}

For the reader's convenience, we outline the proof of this lemma, which follows
the proof of the Hadamard representation of entire functions of finite order
of growth.

\medskip\noindent{\em Proof:}
In order to prove the first statement, we fix $z\ne 0$, write
\[
U_R(z) = \Bigl(\, \int_{|\zeta|\le 2|z|} + \int_{|\zeta| > 2|z|} \, \Bigr)
\mathsf H_s(z/\zeta)\, {\rm d}\mu_R(\zeta)\,,
\]
and estimate each of the integrals starting with the second one.
For $|w|\le \tfrac12$, we have
\[
\mathsf H_s(w) = - {\rm Re\, }\frac{w^{s+1}}{s+1} + O\bigl( |w|^{s+2} \bigr)\,.
\]
Since
\[
\int_{2|z|}^\infty \Bigl( \frac{|z|}t \Bigr)^{s+2} t^{p-1}\, {\rm d}t
= \frac{2^{p-s-2}}{s+2-p}\, |z|^p\,,
\]
the $O$-term presents no trouble. The same can be said about $-{\rm Re\,} w^{s+1}/(s+1)$ if $s+1>p$, that is, if $p$ is non-integer.

If $p$ is a positive integer and $s+1=p$, observing that by~\eqref{eq:Lindelof},
\[
\int_{|\zeta|>2|z|} \zeta^{-p}\, {\rm d} \mu_R (\zeta) = 0, \quad \text{for\ all\ }
0<R<\infty,
\]
we see that the term $-{\rm Re\,} w^{s+1}/(s+1)$ has no influence at all.

It remains to crudely estimate the integrals over $|\zeta|\le 2|z|$:
\begin{align*}
\int_0^{2|z|} \log\Bigl| 1 - \frac{z}t \Bigr|\, t^{p-1}\, {\rm d} t
&\le \int_0^{2|z|} \log \Bigl( 1+\frac{|z|}t \Bigr)\, t^{p-1}\, {\rm d}t  \\
&= |z|^p\, \int_0^2 \log\Bigl( 1+\frac1{t} \Bigr)\,t^{p-1}\, {\rm d}t
\end{align*}
and
\[
\int_0^{2|z|} \Bigl( \frac{|z|}t \Bigr)^j\,t^{p-1}\, {\rm d}t
= |z|^p\, \int_0^2 t^{p-1-j}\, {\rm d}t, \quad 1 \le j \le s<p\,.
\]
Thus, we have our growth bound.

Similar estimates show that $U_{R_2}(z)-U_{R_1}(z)$ tends to zero locally uniformly in $z$, as
$R_2>R_1$ and $R_1\to\infty$, which proves the locally uniform convergence of $U_R(z)$ as $R\to\infty$. The limiting
function $U$ is continuous and, for every $f\in C_0^\infty$,
\[
\int U \Delta f = \lim_{R\to\infty}\, \int U_R \Delta f
= \lim_{R\to\infty}\, \int f\, {\rm d}\mu_R = \int f\, {\rm d}\mu\,,
\]
that is, $\Delta U = \mu$ in the sense of distributions. Thus,
$\Delta (U-K) = 0$ in the sense of distributions and, by the classical Weyl lemma,
the function $U-K$ is harmonic.
Since it is bounded from above by $C|z|^p$ everywhere in $\mathbb C$, by a version
of the Liouville theorem, $U-K=H$, where $H$ is a real part of a polynomial of degree at most $p$.
Since $U$ and $K$ vanish at the origin, $H$ also vanishes at the origin. Since
$H$ is harmonic, vanishes at the origin, and bounded from above by $O(|z|^p)$,
its absolute value is also bounded by $O(|z|^p)$ everywhere
in $\mathbb C$. This information suffices to conclude that $H(z)$ vanishes identically
when $p$ is non-integer, and that $H(z) = {\rm Re}(cz^p)$ when $p$ is a positive integer.
\hfill $\Box$

\subsection*{Construction of the entire function $S$}
Now  we switch from the continuous measure $\mu$ to the discrete measure
$\eta = \sum_{\zeta\in Z} \delta_\zeta$, where $Z$ is a $k$-smooth set of points.
We assume that $0\notin Z$ (otherwise, we construct $S$ for the set $Z\setminus \{0\}$ as below,
and take the function $zS$).

We will again denote $\eta_R = \done_{R\overline{\mathbb D}}\, \eta$ and let $s$ be an integer
such that $s<p\le s+1$.
Set
\[
S_R(z) = \prod_{\zeta\in Z\cap R\overline{\mathbb D}} \mathsf E_s\Bigl( \frac{z}\zeta \Bigr),
\]
Then
\[
\log|S_R(z)| = \int_{\mathbb C} \mathsf H_s \Bigl( \frac{z}\zeta \Bigr)\, {\rm d}\eta_R(\zeta).
\]
The functions $S_R$ are entire and their zero sets coincide with $Z\cap R\overline{\mathbb D}$.
We want to show that they also converge uniformly on compact sets and that $\log |S_R|$ does not differ
from $U_R$ too much.

For non-integer values of $p$, the convergence follows from the estimate
\[
\mathsf E_s (w) = 1 + O(|w|^{s+1}), \quad |w|\le \tfrac12\,,
\]
combined with convergence of the series $ \displaystyle \sum_{\zeta\in Z} |\zeta|^{-t} $
for any $t>p$, in particular, for $t=s+1$.

When $p$ is an integer, we have $s+1=p$,
\[
\mathsf E_s (w) = 1 - \frac{w^p}p + O(|w|^{p+1}), \quad |w|\le \tfrac12\,.
\]
Denote by $n_j$ the counting measure of the portion of $Z$ lying on the ray $\arg(\zeta)=\theta_j$,
and let $\nu_j (r) = \displaystyle n_j (r\overline{\mathbb D} ) -\frac{m_j}{2\pi p}\, r^p$.
Denoting by $ \Omega $ complex-valued constants,
we have, for $R\ge 1$,
\[
\sum_{\zeta\in Z\cap R\overline{\mathbb D}}\, \frac1{\zeta^{p}}
\underset{\eqref{eq:Lindelof}}= \, \sum_{j=1}^n e^{-{\rm i}p\theta_j}
\int_{[1, R]} \frac{{\rm d}\nu_j(r)}{r^{p}}\,  + \Omega\,.
\]
Integrating by parts, and using that, by $k$-smoothness of $Z$,
$|\nu_j (t)| = O(1)$ on $\R_+$,
we see that the RHS equals $\Omega_1 + O(R^{-p})$ as $R\to\infty$, which proves the
existence of the limit of $S_R (z)$.

It remains to show that the difference of the functions $V_R = \log|S_R|$ and $U_R$ is close to
a harmonic polynomial of degree at most $p$.
We have
\[
\bigl( U_R - V_R)(z) =
- \sum_j \int_{[0, R]} \mathsf H_s\bigl( \frac{ze^{-{\rm i}\theta_j}}{t} \bigr)\, {\rm d}\nu_j(t)
\]
We consider the contribution of the positive semi-axis, the contribution of other rays is
treated similarly. Dropping the index $j$, we have
\begin{multline*}
\int_{[0, R]} \mathsf H_s\bigl( \frac{z}t \bigr)\, {\rm d}\nu (t) \\
= \mathsf H_s\Bigl( \frac{z}R \Bigr)\, \nu (R)
- \int_0^R \nu (t)\,
{\rm Re}\, \Bigl( \frac1{t-z} - \frac1t \Bigr)\, {\rm d}t
- \sum_{k=1}^s {\rm Re}(z^k)
\int_0^R \frac{\nu(t)}{t^{k+1}}\, {\rm d}t.
\end{multline*}
The first term on the RHS tends to zero as $R\to\infty$ locally uniformly in $z$ and can be discarded.
Since $\nu (t) = - c t^p$ near the origin and is bounded near infinity, the coefficients
\[
\int_0^\infty \nu (t) t^{-k-1}\, {\rm d}t, \quad 1\le k \le s,
\]
are convergent integrals, i.e., the third term on the RHS converges to the real part of a
polynomial of degree $s$, as $R\to\infty$. It remains to bound
\[
\int_0^R |\nu (t)|\, \Bigl| \frac1{t-z} - \frac1t \Bigr|\, {\rm d}t\,.
\]
Let $\la_1$ be the smallest absolute value of the points from $Z$ lying on the positive ray. Then,
for $|z|\ge 1+\la_1$, the integral we are estimating is bounded by
\[
c\int_0^{\la_1} \frac{t^p}t\, {\rm d} t +
\| \nu \|_\infty\,
\Bigl[ \int_{\la_1}^{2|z|} \frac{{\rm d}t}t
+ \int_1^{2|z|} \frac{{\rm d}t}{|t-z|} + \int_{2|z|}^\infty \frac{|z|}{|t-z|t}\, {\rm d}t \Bigr]\,.
\]
The third integral needs a little care. Let $z=x+{\rm i}y$.
For $x<0<t$, we have $|t-z|^{-1}\le t^{-1}$, and therefore,
$\displaystyle \int_1^{2|z|} \frac{{\rm d}t}{|t-z|} \le \log 2 +\log |z| $.
For $x\ge 0$, we have
\begin{multline*}
\int_1^{2|z|} \frac{{\rm d}t}{|t-z|} \le
\int_1^{\max(x-|y|,1)} \frac{{\rm d}t}{x-t} + \int_{x-|y|}^{x+|y|} \frac{{\rm d}t}{|y|}
+ \int_{x+|y|}^{2|z|} \frac{{\rm d}t}{t-x} \\
\le
2 \int_{|y|}^{2|z|} \frac{{\rm d}s}{s} + 2
= 2\log |z| + 2\log\frac1{|y|} + 2\log 2 +2\,.
\end{multline*}
Hence,
\begin{multline*}
c\int_0^{\la_1} \frac{t^p}t\, {\rm d} t +
\| \nu \|_\infty\,
\Bigl[ \int_{\la_1}^{2|z|} \frac{{\rm d}t}t
+ \int_1^{2|z|} \frac{{\rm d}t}{|t-z|} + \int_{2|z|}^\infty \frac{|z|}{|t-z|t}\, {\rm d}t \Bigr] \\
\le C\Bigl[ 1 + \log |z| + \log\Bigl(1+\frac1{|y|}\Bigr) \Bigr]\,.
\end{multline*}

Correcting the functions $S_R$ by $e^P$, where $P$ is a polynomial of degree at most $p$,
independent of $R$,
we get the limiting entire function $S$ with the zero set $Z$ satisfying
\begin{equation}\label{eq:estimate-S-K}
\bigl| \log |S(z)| - K(z) \bigr|
\le C \Bigl[ \log (2+|z|) + \log\Bigl(1 + \frac1{d(z)} \Bigr) \Bigr]\,,
\end{equation}
where $d(z)$ is the distance from $z$ to the union of our rays.

Now choose some big power $\tau>p$ and consider a point $z$, $|z|\ge 1$,
with $d(z)<|z|^{-\tau}$
(if $d(z)\ge |z|^{-\tau}$, then~\eqref{eq:estimate-S-K} directly yields the needed bounds of $S$).
Without loss of generality, we may assume that the nearest ray to $z$ is the positive
real semi-axis. Consider the disk $D$ centered at $z$ of radius $\rho=2|z|^{-\tau}$. Since
$K(re^{{\rm i}\theta})=k(\theta)r^p$ with a Lipschitz function $k$,
and $\tau>p$, we have
\begin{equation}\label{eq:K}
\sup\{ |K(z') - K(z)|\colon z'\in D\} = O(1), \quad z\in\mathbb C.
\end{equation}
By subharmonicity of $\log |S|$,
\begin{align*}
\log|S(z)| &\le \frac1{2\pi}\, \int_{-\pi}^\pi \log|S(z+\rho e^{{\rm i}\theta})|\, {\rm d}\theta \\
&\le \frac1{2\pi}\, \int_{-\pi}^\pi K(z+\rho e^{{\rm i}\theta})\, {\rm d}\theta \\
& \qquad \qquad \qquad + C \Bigl[
\frac1{2\pi}\, \int_{-\pi}^\pi \log\Bigl( 1 + \frac1{d(z+\rho e^{{\rm i}\theta})} \Bigr)\,
{\rm d}\theta + \log(1+|z|) \Bigr]  \\
&\underset{\eqref{eq:K}}\le K(z) + C\log (1+|z|)\,.
\end{align*}
This gives us estimate (I) is Theorem~L.

Furthermore, if this disk has no zeroes of $S$ in it, then
\begin{align*}
\log|S(z)| &= \frac1{2\pi}\, \int_{-\pi}^\pi
\log|S(z+\rho e^{{\rm i}\theta})|\, {\rm d}\theta \\
&\ge \frac1{2\pi}\, \int_{-\pi}^\pi K(z+\rho e^{{\rm i}\theta})\, {\rm d}\theta \\
& \qquad \qquad \qquad - C \Bigl[
\frac1{2\pi}\, \int_{-\pi}^\pi \log\Bigl( 1 + \frac1{d(z+\rho e^{{\rm i}\theta})} \Bigr)
\, {\rm d}\theta + \log(1+|z|) \Bigr] \\
&\ge K(z) - C\log (1+|z|)
\end{align*}
holds as well, which yields estimate (II). \hfill $\Box$

\bigskip
\medskip

\noindent Aleksei Kulikov: School of Mathematical Sciences, Tel Aviv University, Tel Aviv, 69978, Israel
\newline {\tt alekseik@tauex.tau.ac.il, lyosha.kulikov@mail.ru}

\medskip\noindent Fedor Nazarov: Department of Mathematical Sciences, Kent State University, Kent, OH 44242, USA \newline {\tt nazarov@math.kent.edu}

\medskip\noindent Mikhail Sodin:
School of Mathematical Sciences, Tel Aviv University, Tel Aviv, 69978, Israel
\newline {\tt sodin@tauex.tau.ac.il}

\end{document}